\documentclass{article}

\usepackage[final, nonatbib]{neurips_2022}

\usepackage[utf8]{inputenc} \usepackage[T1]{fontenc}    \usepackage{hyperref}       \usepackage{url}            \usepackage{booktabs}       \usepackage{amsfonts}       \usepackage{nicefrac}       \usepackage{microtype}      \usepackage{xcolor}         
\usepackage{tikz}
\usepackage{nicematrix}
\usetikzlibrary{positioning}
\usepackage{bm}
\usepackage{xcolor,colortbl}
\usepackage{lmodern}
\usepackage{scalerel}
\usepackage{tablefootnote}
\usepackage{microtype}
\usepackage{graphicx}
\usepackage{subfigure}
\usepackage{booktabs}
\usepackage{hyperref}
\usepackage{cleveref}
\usepackage{amssymb}
\usepackage{algorithm,algorithmic}
\usepackage{wrapfig}

\newtheorem{theorem}{Theorem}
\newtheorem{definition}{Definition}
\newtheorem{proposition}{Proposition}
\newtheorem{corollary}{Corollary}

\newtheorem{remark}{Remark}

\newcommand{\X}{\bm{X}}
\newcommand{\Z}{\bm{Z}}
\newcommand{\W}{\bm{W}}

\newcommand{\ind}{1\!\!1}

\title{A Probabilistic Graph Coupling View of Dimension Reduction}

\author{
  Hugues Van Assel \\
  UMPA, CNRS \\
  ENS Lyon\\
  \texttt{hugues.van\_assel@ens-lyon.fr}
  \And
  Thibault Espinasse \\
  Institut Camille Jordan Lyon 1 \\
  Inria Dracula \\
  \texttt{espinasse@math.univ-lyon1.fr}
  \And
  Julien Chiquet \\
  AgroParisTech \\
  INRAE \\
  \texttt{julien.chiquet@inrae.fr} \\
  \And
  Franck Picard \\
  LBMC, CNRS \\
  ENS Lyon \\
  \texttt{franck.picard@ens-lyon.fr}
}

\begin{document}

\maketitle

\begin{abstract}
  Most popular dimension reduction (DR) methods like t-SNE and UMAP are based on minimizing a cost between input and latent pairwise similarities. Though widely used, these approaches lack clear probabilistic foundations to enable a full understanding of their properties and limitations. To that extent, we introduce a unifying statistical framework based on the coupling of hidden graphs using cross-entropy. These graphs induce a Markov random field dependency structure among the observations in both input and latent spaces. We show that existing pairwise similarity DR methods can be retrieved from our framework with particular choices of priors for the graphs. Moreover, this reveals that these methods relying on shift-invariant kernels suffer from a statistical degeneracy that explains poor performances in conserving coarse-grain dependencies. New links are drawn with PCA which appears as a non-degenerate graph coupling model.
\end{abstract}

\section{Introduction}\label{intro}

Dimensionality reduction (DR) is of central importance when dealing with high-dimensional data \cite{donoho2000high}. It mitigates the curse of dimensionality, allowing for greater statistical flexibility and less computational complexity. DR also enables visualization that can be of great practical interest for understanding and interpreting the structure of large datasets.
Most seminal approaches include Principal Component Analysis (PCA) \cite{pearson1901liii},  multidimensional scaling \cite{kruskal1978multidimensional} and more broadly kernel eigenmaps methods such as Isomap \cite{balasubramanian2002isomap}, Laplacian eigenmaps \cite{belkin2003laplacian} and diffusion maps \cite{coifman2006diffusion}. These methods share the definition of a pairwise similarity kernel that assigns a high value to close neighbors and the resolution of a spectral problem. They are well understood and unified in the kernel PCA framework \cite{ham2004kernel}.

In the past decade, the field has witnessed a major shift with the emergence of a new class of methods. They are also based on pairwise similarities but these are not converted into inner products. Instead, they define pairwise similarity functions in both input and latent spaces and optimize a cost between the two. Among such methods, the Stochastic Neighbor Embedding (SNE) algorithm \cite{NIPS2002SNE}, its heavy-tailed symmetrized version t-SNE \cite{maaten2008tSNE} or more recent approaches like LargeVis \cite{tang2016visualizing} and UMAP \cite{mcinnes2018umap} are arguably the most used in practice. These will be referred to as \textit{SNE-like} or \textit{neighbor embedding} methods in what follows. They are increasingly popular and now considered as the state-of-art techniques in many fields \cite{li2017application,kobak2019art,anders2018dissecting}. Their popularity is mainly due to their exceptional ability to preserve local structure, \textit{i.e.}\ close points in the input space have close embeddings, as shown empirically \cite{wang2021understanding}. They also demonstrate impressive performances in identifying clusters \cite{arora2018analysis, linderman2019clustering}. However this is done at the expense of global structure, that these methods struggle in preserving \cite{wattenberg2016use, coenen2019understanding} \textit{i.e.}\ the relative large-scale distances between embedded points do not necessarily correspond to the original ones. 

Due to a lack of clear probabilistic foundations, these properties remain mostly empirical. This gap between theory and practice is detrimental as practitioners may rely on strategies that are not optimal for their use case.
While recent software developments are making these methods more scalable \cite{chan2018t,pezzotti2019gpgpu,linderman2019fast} and further expanding their use, the need for a well-established probabilistic framework is becoming more prominent.
In this work, we define the generative probabilistic model that encompasses current embedding methods while establishing new links with the well-established PCA model.

\paragraph{Outline.} 
Consider $\X = (\X_1,..., \X_n)^\top \in \mathbb{R}^{n \times p}$, an input dataset that consists of $n$ vectors of dimension $p$. Our task is to embed $\X$ in a lower dimensional space of dimension $q<p$ (typically $q=2$ for visualization), and we denote by $\Z = (\Z_1, ..., \Z_n)^\top \in \mathbb{R}^{n \times q}$ the unknown embeddings. The rationale of our framework is to suppose that the observations $\X$ and $\Z$ are structured by two latent graphs with $\W_{\scaleto{X}{4pt}}$ and $\W_{\scaleto{Z}{4pt}}$ standing for their $n$-square weight matrices.
As the goal of DR is to preserve the input's structure in the latent space, we propose to find the best low-dimensional representation $\Z$ of $\X$ such that $\W_{\scaleto{X}{4pt}}$ and $\W_{\scaleto{Z}{4pt}}$ are close. To build a flexible and robust probabilistic framework, we consider random graphs distributed according to some predefined prior distributions. Our objective is to match the posterior distributions of $\W_{\scaleto{X}{4pt}}$ and $\W_{\scaleto{Z}{4pt}}$. Note that as they share the same dimensionality the latter graphs can be easily compared unlike $\X$ and $\Z$. The coupling is done with a cross-entropy criterion, the minimization of which will be referred to as graph coupling.

In this work, our main contributions are as follows.

\begin{itemize}
    \item We show that SNE, t-SNE, LargeVis and UMAP are all instances of graph coupling and characterized by different choices of prior for discrete latent structuring graphs (\cref{sec:GC_unified}). We demonstrate that such graphs essentially capture conditional independencies among rows through a pairwise Markov Random Field (MRF) model whose construction can be found in \cref{sec:graph_structure}.
    \item We uncover the intrinsic probabilistic property explaining why such methods perform poorly on conserving the large-scale structure of the data as a consequence of a degeneracy of the MRF when shift-invariant kernels are used (\cref{prop:integrability_pairwise_MRF}). Such degeneracy induces the loss of the relative positions of clusters corresponding to the connected components of the posterior latent graphs which distributions are identified (\cref{prop:posterior_W}). These findings are highlighted by a new initialization of the embeddings that is empirically tested (\cref{sec:towards_large_scale}).
    \item We show that for Gaussian MRFs, when adapting graph coupling to precision matrices with suitable priors, PCA appears as a natural extension of the coupling problem in its continuous version (\cref{PCA_graph_coupling}). Such model does not suffer from the aforementioned degeneracy and hence preserves the large-scale structure.
\end{itemize}

\section {Shift-Invariant Pairwise MRF to Model Row Dependencies} \label{sec:graph_structure}

We start by defining the distribution of the observations given a graph. The latter takes the form of a pairwise MRF model which as we show is improper (\textit{i.e.}\ not integrable on $\mathbb{R}^{n \times p}$) when shift-invariant kernels are used. We consider a fixed directed graph $\W \in \mathcal{S}_{\scaleto{W}{4pt}}$ where:
$$\mathcal{S}_{\scaleto{W}{4pt}} = \left\{\W \in \mathbb{N}^{n \times n} \mid \forall (i,j) \in [n]^2, W_{ii}=0, W_{ij} \leq n \right\}$$
Throughout, $(E, \mathcal{B}(E), \lambda_E)$ denotes a measure space where $\mathcal{B}(E)$ is the Borel $\sigma$-algebra on $E$ and $\lambda_E$ is the Lebesgue measure on $E$.

\subsection{Graph Laplacian Null Space}\label{sec:laplacian_prop}
A central element in our construction is the graph Laplacian linear map, defined as follows, where $\mathcal{S}^n_+(\mathbb{R})$ is the set of positive semidefinite matrices.
\begin{definition}\label{graph_laplacian}
The graph Laplacian operator is the map $L \colon \mathbb{R}_+^{n \times n} \to \mathcal{S}^n_+(\mathbb{R})$ such that
$$\text{for } (i,j) \in [n]^2, \quad L(\W)_{ij} = \left\{
\begin{array}{ll}
    - W_{ij} & \text{if } i \neq j \\
    \sum_{k \in [n]} W_{ik} & \text{otherwise} \:.
\end{array} 
\right. $$
\end{definition}
With an abuse of notation, let $\bm{L} = L(\overline{\W})$ where $\overline{\W} = \W + \W^\top$. Let $(C_1,...,C_{\scaleto{R}{4pt}})$ be a partition of $[n]$ (\textit{i.e.}\ the set $\{1,2,...,n\}$) corresponding to the connected components (CCs) of $\overline{\W}$. As well known in spectral graph theory \cite{Chung97}, the null space of $\bm{L}$ is spanned by the orthonormal vectors $\{\bm{U}_{r}\}_{r \in [R]}$ such that for $r \in [R]$,
$\bm{U}_{r} = \left(n_r^{-1/2} \ind_{i \in C_r}\right)_{i \in [n]}$ with $n_r = \operatorname{Card}(C_r)$. By the spectral theorem, $\bm{U}_{\scaleto{[R]}{5pt}}$ can be completed such that $\bm{L} = \bm{U \Lambda U^\top}$ where $\bm{U} = (\bm{U}_1, ..., \bm{U}_n)$ is orthogonal and $\bm{\Lambda} = \operatorname{diag}((\lambda_i)_{i \in [n]})$ with $0 = \lambda_1 = ... = \lambda_R < \lambda_{R+1} \leq ... \leq \lambda_n$. In the following, the data is split into two parts: $\X_{\scaleto{M}{4pt}}$, the orthogonal projection of $\X$ on $\mathcal{S}_{\scaleto{M}{4pt}} = (\ker \bm{L}) \otimes \mathbb{R}^p$, and $\X_{\scaleto{C}{4pt}}$, the projection on $\mathcal{S}_{\scaleto{C}{4pt}} = (\ker \bm{L})^{\perp} \otimes \mathbb{R}^p$. For $i \in [n]$, $\X_{\scaleto{M}{4pt},i} = \sum_{r \in [R]} n_r^{-1} \ind_{i \in C_r}\sum_{\ell \in C_r} \X_{\ell} $ hence $\X_{\scaleto{M}{4pt}}$ stands for the empirical means of $\X$ on CCs, thus modelling the CC positions, while $\X_{\scaleto{C}{4pt}} = \X - \X_{\scaleto{M}{4pt}}$ is CC-wise centered, thus modeling the relative positions of the nodes within CCs. We now introduce the probability distribution of these variables.

\subsection{Pairwise MRF and Shift-Invariances}\label{sec:within_CC}

In this work, the dependency structure among rows of the data is governed by a graph. The strength of the connection between two nodes is given by a symmetric function $k : \mathbb{R}^p \to \mathbb{R}_+$. We consider the following pairwise MRF unnormalized density function:
\begin{align}\label{eq:unnormalized_MRF}
  f_{k} \colon (\X,\W) &\mapsto \prod_{(i,j) \in [n]^2} k(\X_{i} - \X_{j})^{W_{ij}} \: .
\end{align}
As we will see shortly, the above is at the heart of DR methods based on pairwise similarities. Note that as $k$ measures the similarity between couples of samples, $f_k$ will take high values if the rows of $\X$ vary smoothly on the graph $\W$. Thus we can expect $\X_i$ and $\X_j$ to be close if there is an edge between node $i$ and node $j$ in $\W$. A key remark is that $f_{k}$ is kept invariant by translating $\X_{\scaleto{M}{4pt}}$. Namely for all $\X \in \mathbb{R}^{n \times p}$, $f_{k}(\X, \W) = f_{k}(\X_{\scaleto{C}{4pt}}, \W)$. This invariance results in $f_{k}(\cdot, \W)$ being non-integrable on $\mathbb{R}^{n \times p}$, as we see with the following example. 

\paragraph{Gaussian kernel.} For a positive definite matrix $\bm{\Sigma} \in \mathcal{S}^n_{++}(\mathbb{R})$, consider the Gaussian kernel $k : \bm{x} \mapsto e^{- \frac{1}{2}\| \bm{x} \|_{\bm{\Sigma}}^2}$ where $\bm{\Sigma}$ stands for the covariance among columns. One has:
\begin{align}\label{eq:gaussian_kernel}
    \log f_{k}(\X, \W) &= -\sum_{(i,j) \in [n]^2} W_{ij} \| \bm{X}_{i}-\bm{X}_{j} \|^2_{\bm{\Sigma}}
    = - \operatorname{tr} \left(\bm{\Sigma}^{-1} \bm{X}^{T} \bm{L} \bm{X}\right)
\end{align}
by property of the graph Laplacian (\cref{graph_laplacian}). In this case, it is clear that due to the rank deficiency of $\bm{L}$, $f_{k}(\cdot, \W)$ is only $\lambda_{\mathcal{S}_{\scaleto{C}{3pt}}}$-integrable. In general DR settings, one does not want to rely on Gaussian kernels only. A striking example is the use of the Student kernel in t-SNE \cite{maaten2008tSNE}. Heavy-tailed kernels appear useful when the dimension of the embeddings is smaller than the intrinsic dimension of the data \cite{kobak2019heavy}. Our contribution provides flexibility by extending the previous result to a large class of kernels, as stated in the following theorem.

\begin{theorem}\label{prop:integrability_pairwise_MRF}
If $k$ is $\lambda_{\mathbb{R}^p}$-integrable and bounded above $\lambda_{\mathbb{R}^p}$-almost everywhere then $f_{k}(\cdot, \W)$ is $\lambda_{\mathcal{S}_{\scaleto{C}{3pt}}}$-integrable.
\end{theorem}

We refer to \cref{proof:lambda_perp_integrability} for the proof.
We can now define a distribution on $(\mathcal{S}_{\scaleto{C}{4pt}}, \mathcal{B}(\mathcal{S}_{\scaleto{C}{4pt}}))$, where $\mathcal{C}_{k}(\W) = \int f_{k}(\cdot, \W) d\lambda_{\mathcal{S}_{\scaleto{C}{3pt}}}$:
\begin{align}\label{eq:proba_perp}
\mathbb{P}_{k}(d\X_{\scaleto{C}{4pt}} | \W) = \mathcal{C}_{k}(\W)^{-1} f_{k}(\X_{\scaleto{C}{4pt}}, \W) \lambda_{\mathcal{S}_{\scaleto{C}{3pt}}}(d\X_{\scaleto{C}{4pt}}) \: .
\end{align}

\begin{remark}
Kernels may have node-specific bandwidths $\bm{\tau}$, set during a pre-processing step, giving $f_{k}(\X,\W) = \prod_{(i,j)} k((\X_{i} - \X_{j})/\tau_{i})^{W_{ij}}$. Note that such bandwidth does not affect the degeneracy of the distribution and \cref{prop:integrability_pairwise_MRF} still holds.
\end{remark}

\paragraph{Between-Rows Dependency Structure.} By symmetry of $k$, reindexing gives: $f_{k}(\X, \W) = \prod_{j \in [n]} \prod_{i \in [j]} k(\X_{i} - \X_{j})^{\overline{W}_{ij}}$. Hence distribution \eqref{eq:proba_perp} boils down to a pairwise MRF model \cite{clifford1990markov} with respect to the undirected graph $\overline{\W}$, $\mathcal{C}_{k}$ playing the role of the partition function. Note that since $f_k$ (Equation \ref{eq:unnormalized_MRF}) trivially factorize according to the cliques of $\overline{\W}$, the Hammersley-Clifford theorem ensures that the rows of $\X_{\scaleto{C}{4pt}}$ satisfy the local and global Markov properties with respect to $\overline{\W}$. 

\subsection{Uninformative Model for CC-wise Means}

We showed that the MRF (\ref{eq:unnormalized_MRF}) is only integrable on $\mathcal{S}_{\scaleto{C}{4pt}}$, the definition of which depends on the connectivity structure of $\W$. As we now demonstrate, the latter MRF can be seen as a limit of proper distributions on $\mathbb{R}^{n \times p}$, see \textit{e.g.}\ \cite{rue2005gaussian} for a similar construction in the Gaussian case. 
We introduce the Borel function $f^{\varepsilon}(\cdot, \W) \colon \mathbb{R}^{n \times p} \to \mathbb{R}_+$ for $\varepsilon > 0$ such that for all $\X \in \mathbb{R}^{n \times p}$, $f^{\varepsilon}(\X, \W) = f^{\varepsilon}(\X_{\scaleto{M}{4pt}}, \W)$. To allow $f^{\varepsilon}$ to become arbitrarily non-informative, we assume that for all $\W \in \mathcal{S}_{\scaleto{W}{4pt}}$, $f^\varepsilon(\cdot, \W)$ is $\lambda_{\scaleto{\mathcal{S}_{\scaleto{M}{3pt}}}{6pt}}$-integrable for all $\varepsilon \in \mathbb{R}^*_+$ and $f^{\varepsilon}(\cdot, \W) \xrightarrow[\varepsilon \to 0]{} 1$ almost everywhere.
We now define the conditional distribution on $(\mathcal{S}_{\scaleto{M}{4pt}}, \mathcal{B}(\mathcal{S}_{\scaleto{M}{4pt}}))$ as follows:
\begin{align}\label{eq:proba_parallel}
     \mathbb{P}^{\varepsilon}(d\X_{\scaleto{M}{4pt}}| \W) = \mathcal{C}^{\varepsilon}(\W)^{-1} f^{\varepsilon}(\X_{\scaleto{M}{4pt}}, \W) \lambda_{\mathcal{S}_{\scaleto{M}{3pt}}}(d\X_{\scaleto{M}{4pt}})
\end{align}
where $\mathcal{C}^{\varepsilon}(\W) = \int f^{\varepsilon}(\cdot, \W) d\lambda_{\mathcal{S}_{\scaleto{M}{3pt}}}$.
With this at hand, the joint conditional is defined as the product measure of (\ref{eq:proba_perp}) and (\ref{eq:proba_parallel}) over the row axis, the integrability of which is ensured by the Fubini-Tonelli theorem. In the following we will use the compact notation $\mathcal{C}^{\varepsilon}_k(\W) = \mathcal{C}_k(\W)\mathcal{C}^{\varepsilon}(\W)$ for the joint normalizing constant.

\begin{remark}
At the limit $\varepsilon \to 0$ the above construction amounts to setting an infinite variance on the distribution of the empirical means of $\X$ on CCs, thus losing the inter-CC structure. 
\end{remark}

As an illustration, one can structure the CCs' relative positions according to a Gaussian model with positive definite precision $\varepsilon \bm{\Theta} \in \mathcal{S}_{++}^R(\mathbb{R})$, as it amounts to choosing $f^{\varepsilon} : \X \to \exp \left(-\frac{\varepsilon}{2} \operatorname{tr}\left(\bm{\Sigma}^{-1}\X^\top\bm{U}_{\scaleto{[:R]}{5pt}}  \bm{\Theta}\bm{U}^\top_{\scaleto{[R]}{5pt}} \X\right)\right)$ such that: $\mathrm{vec}(\X_{\scaleto{M}{4pt}}) | \bm{\Theta} \sim \mathcal{N}\left(\bm{0}, \left(\varepsilon \bm{U}_{\scaleto{[:R]}{5pt}}  \bm{\Theta}\bm{U}^\top_{\scaleto{[R]}{5pt}}\right)^{-1} \otimes \bm{\Sigma}\right)$ where $\otimes$ denotes the Kronecker product.

\section{Graph Coupling as a Unified Objective for Pairwise Similarity Methods}\label{sec:GC_unified}

In this section, we show that neighbor embedding methods can be recovered in the presented framework. They are obtained, for particular choices of graph priors, at the limit $\varepsilon \to 0$ when $f^{\varepsilon}$ becomes non-informative and the CCs' relative positions are lost. 

We now turn to the priors for $\W$. Our methodology is similar to that of constructing conjugate priors for distributions in the exponential family \cite{wainwright2008graphical}, notably we insert the cumulant function $\mathcal{C}_k^{\varepsilon}$ (\textit{i.e.}\ normalizing constant of the conditional) as a multivariate term of the prior. 
We consider different forms: binary ($B$), unitary out-degree ($D$) and $n$-edges ($E$), relying on an additional term ($\Omega$) to constrain the topology of the graph. For a matrix $\bm{A}$, $A_{i+}$ denotes $\sum_j A_{ij}$ and $A_{++}$ denotes $\sum_{ij} A_{ij}$. In the following, $\bm{\pi}$ plays the role of the edge's prior. The latter can be leveraged to incorporate some additional information about the dependency structure, for instance when a network is observed \cite{li2020high}. 

\begin{definition}\label{def:prior_W}
Let $\bm{\pi} \in \mathbb{R}_+^{n \times n}$, $\varepsilon \in \mathbb{R}_+$, $\alpha \in \mathbb{R}$, $k$ satisfies the assumptions of \cref{prop:integrability_pairwise_MRF} and $\mathcal{P} \in \{B,D,E\}$. For $\W \in \mathcal{S}_{\scaleto{W}{4pt}}$ we introduce:
$$\mathbb{P}_{\scaleto{\mathcal{P},k}{5pt}}^{\varepsilon}(\bm{W}; \bm{\pi}, \alpha) \propto \mathcal{C}^{\varepsilon}_k(\W)^{\alpha} \: \Omega_{\scaleto{\mathcal{P}}{4pt}}(\W) \prod_{(i,j) \in [n]^2} \pi_{ij}^{W_{ij}}$$
where $\Omega_{\scaleto{B}{4pt}}(\W) = \prod_{ij} \ind_{W_{ij} \leq 1}$, $\Omega_{\scaleto{D}{4pt}}(\W) = \prod_{i} \ind_{W_{i+} = 1}$ and $\Omega_{\scaleto{E}{4pt}}(\W) = \ind_{W_{++} = n}\prod_{ij}(W_{ij}!)^{-1}$.
\end{definition}

When $\alpha = 0$, the above no longer depends on $\varepsilon$ and $k$. We will use the compact notation $\mathbb{P}_{\scaleto{\mathcal{P}}{4pt}}(\W ; \bm{\pi}) = \mathbb{P}_{\scaleto{\mathcal{P},k}{5pt}}^{\varepsilon}(\bm{W}; \bm{\pi}, 0)$. Note that by $\W \sim \mathbb{P}_{\scaleto{\mathcal{P}}{4pt}}(\cdot \: ; \bm{\pi})$ we have the following simple Bernoulli $(\mathcal{B})$ and multinomial $(\mathcal{M})$ distributions, where matrix or vector division is to be understood as element-wise.
\begin{itemize}
    \item If $\mathcal{P} = B$, $\forall (i,j) \in [n]^2, \: W_{ij} \stackrel{\perp\!\!\!\!\perp}{\sim} \mathcal{B}\left(\pi_{ij}/(1 + \pi_{ij}) \right)$.
    \item If $\mathcal{P} = D$, $\forall i \in [n], \: \W_{i} \stackrel{\perp\!\!\!\!\perp}{\sim} \mathcal{M}\left(1, \bm{\pi}_{i}/\pi_{i+} \right)$.
    \item If $\mathcal{P} = E$, $\W \sim \mathcal{M}\left(n, \bm{\pi}/\pi_{++} \right)$.
\end{itemize}

We now show that the posterior distribution of the graph given the observations takes a simple form when the distribution of CC empirical means $\bm{X}_{\scaleto{M}{4pt}}$ diffuses \textit{i.e.}\ when $\varepsilon \to 0$ (a proof of the following result can be found in \cref{proof:posterior_limit}). In the following, $\odot$ stands for the Hadamard product and $\mathcal{D}$ for the convergence in distribution.

\begin{proposition}\label{prop:posterior_W}
Let $\bm{\pi} \in \mathbb{R}_+^{n \times n}$, $k$ satisfies the assumptions of \cref{prop:integrability_pairwise_MRF} with  $\bm{K}_{\scaleto{X}{4pt}} = (k(\X_{i} - \X_{j}))_{(i,j) \in [n]^2}$ and $\mathcal{P}\in \{B, D, E\}$. If $\W^{\varepsilon} \sim \mathbb{P}_{\scaleto{\mathcal{P},k}{5pt}}^{\varepsilon}(\cdot \: ; \bm{\pi},1)$ then
$$\W^{\varepsilon} | \X \xrightarrow[\varepsilon \to 0]{\mathcal{D}} \mathbb{P}_{\scaleto{\mathcal{P}}{4pt}}(\cdot \: ;\bm{\pi} \odot \bm{K}_{\scaleto{X}{4pt}}) \:.$$
\end{proposition}

\begin{remark}
For all $\W \in \mathcal{S}_{\scaleto{W}{4pt}}$, $\mathcal{C}^{\varepsilon}(\W)$ diverges as $\varepsilon \to 0$, hence the graph prior (\cref{def:prior_W}) is improper at the limit. This compensates for the uninformative diffuse conditional and allows to retrieve a well-defined tractable posterior limit.
\end{remark}

\subsection{Retrieving Well Known DR Methods}\label{sec:retrieving_DR_methods}

We now provide a unified view of neighbor embedding objectives as a coupling between graph posterior distributions. To that extent we derive the cross entropy associated with the various graph priors at hand. In what follows, $k_x$ and $k_z$ satisfy the assumptions of \cref{prop:integrability_pairwise_MRF} and we denote by $\bm{K}_{\scaleto{X}{4pt}}$ and $\bm{K}_{\scaleto{Z}{4pt}}$ the associated kernel matrices on  $\X$ and $\Z$ respectively. For both graph priors we consider the parameters $\bm{\pi}=\bm{1}$ and $\alpha=1$. For $(\mathcal{P}_{\scaleto{X}{4pt}}, \mathcal{P}_{\scaleto{Z}{4pt}}) \in \{B,D,E\}^2$, we introduce the 
cross entropy between the limit posteriors at $\varepsilon \to 0$,
\begin{align*}
    \mathcal{H}_{\scaleto{\mathcal{P}_X}{4pt}, \scaleto{\mathcal{P}_Z}{4pt}} = - \mathbb{E}_{\W_{\scaleto{X}{3pt}} \sim \mathbb{P}_{\scaleto{\mathcal{P}_{\scaleto{X}{2pt}}}{3pt}}(\cdot;\bm{K}_{\scaleto{X}{3pt}})}[\log \mathbb{P}_{\scaleto{\mathcal{P}_{\scaleto{Z}{3pt}}}{4pt}}(\W_{\scaleto{Z}{4pt}} = \W_{\scaleto{X}{4pt}}; \bm{K}_{\scaleto{Z}{4pt}})]
\end{align*}
defining a coupling criterion to be optimized with respect to embedding coordinates $\Z$. We now go through each couple $(\mathcal{P}_{\scaleto{X}{4pt}}, \mathcal{P}_{\scaleto{Z}{4pt}})$ such that $\operatorname{supp}\left(\mathbb{P}_{\scaleto{\mathcal{P}_X}{4pt}}\right) \subset \operatorname{supp}\left(\mathbb{P}_{\scaleto{\mathcal{P}_Z}{4pt}}\right)$ for the cross-entropy to be defined.

\paragraph{SNE.}
When $\mathcal{P}_{\scaleto{X}{4pt}} = \mathcal{P}_{\scaleto{Z}{4pt}} = D$, the probability of the limit posterior graphs factorizes over the nodes and the cross-entropy between limit posteriors takes the form of the objective of SNE \cite{hinton2002stochastic}, where for $i \in [n], \bm{P}^{D}_{i} = \bm{K}_{\scaleto{X}{4pt},i} / K_{\scaleto{X}{4pt},i+}$ and $\bm{Q}^{D}_{i} = \bm{K}_{\scaleto{Z}{4pt},i} / K_{\scaleto{Z}{4pt},i+}$,
$$\mathcal{H}_{D,D}= - \sum_{i \neq j} P^{D}_{ij} \log Q^{D}_{ij} \:.$$

\paragraph{Symmetric-SNE.}
Choosing $\mathcal{P}_{\scaleto{X}{4pt}} = D$ and $\mathcal{P}_{\scaleto{Z}{4pt}} = E$, we define for $(i,j) \in [n]^2$, $\bm{Q}^{E}_{ij} = K_{\scaleto{Z}{4pt},ij} / K_{\scaleto{Z}{4pt},++}$ and $\overline{P}^{D}_{ij} = P^{D}_{ij} + P^{D}_{ji}$. The symmetry of $\bm{Q}^{E}$ yields:
\begin{align*}
    \mathcal{H}_{D,E} &= - \sum_{i \neq j} P^{D}_{ij} \log Q^{E}_{ij} = - \sum_{i < j} \overline{P}^{D}_{ij} \log Q^{E}_{ij}
\end{align*}
and the symmetrized objective of t-SNE \cite{maaten2008tSNE} is recovered. 

\paragraph{LargeVis.}
Now choosing $\mathcal{P}_{\scaleto{X}{4pt}} = D$ and $\mathcal{P}_{\scaleto{Z}{4pt}} = B$, one can also notice that $\bm{Q}^{B} = \left( K_{\scaleto{Z}{4pt},ij} / (1+K_{\scaleto{Z}{4pt},ij}) \right)_{(i,j) \in [n]^2}$ is symmetric. With this at hand the limit cross-entropy reads
\begin{align*}
    \mathcal{H}_{D,B} = - \sum_{i \neq j} P^{D}_{ij} \log Q^{B}_{ij} + \left(1 - P^{D}_{ij} \right) \log\left(1-Q^{B}_{ij} \right)
    = - \sum_{i < j} \overline{P}^{D}_{ij} \log Q^{B}_{ij} + \left(2-\overline{P}^{D}_{ij}\right) \log (1- Q^{B}_{ij})
\end{align*}
which is the objective of LargeVis \cite{tang2016visualizing}.

\paragraph{UMAP.}
Let us take $\mathcal{P}_{\scaleto{X}{4pt}} = \mathcal{P}_{\scaleto{Z}{4pt}} = B$ and consider the symmetric thresholded graph $\widetilde{\W}_{\scaleto{X}{4pt}} = \ind_{\W_{\scaleto{X}{4pt}} + \W_X^\top \geq 1}$. By independence of the edges, $\widetilde{W}_{\scaleto{X}{4pt},ij} \sim \mathcal{B}\left(\widetilde{P}^{B}_{ij} \right) \quad \text{where} \quad  \widetilde{P}^{B}_{ij} = P^{B}_{ij} + P^{B}_{ji} - P^{B}_{ij} P^{B}_{ji}$ and $\bm{P}^{B} = \left( K_{\scaleto{X}{4pt},ij} / (1+K_{\scaleto{X}{4pt},ij}) \right)_{(i,j) \in [n]^2}$. Coupling $\widetilde{\W}_{\scaleto{X}{4pt}}$ and $\W_{\scaleto{Z}{4pt}}$ gives:
\begin{align*}
    \mathcal{H}_{\widetilde{B},B} &= -2 \sum_{i<j} \widetilde{P}_{ij}^{B} \log Q_{ij}^{B} + \left(1 - \widetilde{P}_{ij}^{B} \right) \log \left( 1 - Q_{ij}^{B} \right)
\end{align*}
which is the loss function considered in UMAP \cite{mcinnes2018umap}, the construction of $\widetilde{\W}_{\scaleto{X}{4pt}}$ being borrowed from section 3.1 of the paper.

\begin{remark}
One can also consider $\mathcal{H}_{E,E}$ but as detailed in \cite{maaten2008tSNE}, this criterion fails at positioning outliers and is therefore not considered. 
Interestingly, any other feasible combination of the presented priors relates to an existing method.
\end{remark}

\subsection{Interpretations}\label{sec:interpretations}

\begin{table}[]
    \caption{Prior distributions for $\W_{\scaleto{X}{4pt}}$ and $\W_{\scaleto{Z}{4pt}}$ associated with the pairwise similarity coupling DR algorithms. Grey-colored boxes are such that the cross-entropy is undefined.}
    \begin{center}
    \begin{small}
    \begin{sc}
    \centering
    \renewcommand{\arraystretch}{2}
    \begin{NiceTabular}{|W{c}{1cm}|W{c}{1.8cm}|W{c}{1.8cm}|W{c}{1.8cm}|}
    \hline
    \diagbox{{\fontsize{12}{15}\selectfont $\mathcal{P}_{\scaleto{X}{4pt}}$}}{{\fontsize{12}{15}\selectfont $\mathcal{P}_{\scaleto{Z}{4pt}}$}} & $B$ & $D$ & $E$ \\
    \hline
    $\widetilde{B}$ & UMAP & \cellcolor{black!10} & \cellcolor{black!10} \\
    \hline
    $D$ & LargeVis & SNE & T-SNE\\
    \hline
    \end{NiceTabular}
    \label{tableau_priors}
    \end{sc}
    \end{small}
    \end{center}
    \label{priors_methods}
\end{table}

As we have seen in \cref{sec:retrieving_DR_methods}, SNE-like methods can all be derived from the graph coupling framework.  What characterizes each of them is the choice of priors considered for the latent structuring graphs. To the best of our knowledge, the presented framework is the first that manages to unify all these DR algorithms. Such a framework opens many perspectives for improving upon current practices as we discuss in \cref{sec:towards_large_scale} and \cref{Perspectives}. 
We now focus on a few insights that our work provides about the empirical performances of these methods. 

\paragraph{Repulsion \& Attraction.}
Decomposing $\mathcal{H}_{\scaleto{\mathcal{P}_X}{4pt}, \scaleto{\mathcal{P}_Z}{4pt}}$ with Bayes' rule and simplifying constant terms one has the following optimization problem: 
\begin{align}\label{eq:optim_H_Z}
    \min_{\Z \in \mathbb{R}^{n \times q}} -\hspace{-0.2cm}\sum_{(i,j) \in [n]^2} \bm{P}^{\scaleto{\mathcal{P}_X}{4pt}}_{ij}\log k_z(\Z_i - \Z_j) + \log \mathbb{P}(\Z).
\end{align}
The first and second terms in \cref{eq:optim_H_Z} respectively summarize the attractive and repulsive forces of the objective. Recall from \cref{prop:posterior_W}
that $\bm{P}^{\scaleto{\mathcal{P}_X}{4pt}}$ is the posterior expectation of $\W_{\scaleto{X}{4pt}}$. Hence in SNE-like methods, the attractive forces resume to a pairwise MRF log likelihood with respect to a graph posterior expectation given $\X$. For instance if $k_z$ is the Gaussian kernel, this attractive term reads $\operatorname{tr} \left(\Z^\top \bm{L}^\star \Z \right)$ where $\bm{L}^\star = \mathbb{E}_{\W \sim \mathbb{P}_{\scaleto{\mathcal{P}_X}{4pt}}(\cdot;\bm{K}_{\scaleto{X}{4pt}})}[L(\W)]$, boiling down to the objective of Laplacian eigenmaps \cite{belkin2003laplacian}. Therefore, for Gaussian MRFs, the attractive forces resume to an unconstrained Laplacian eigenmaps objective. Such link, already noted in \cite{carreira2010elastic}, is easily unveiled in our framework. Moreover, one can notice that only this attractive term depends on $\X$ as the repulsion is given by the marginal term in (\ref{eq:optim_H_Z}). The latter reads $\mathbb{P}(\Z) = \sum_{\W \in \mathcal{S}_{\scaleto{W}{4pt}}} \mathbb{P}(\Z, \W)$ with $\mathbb{P}(\Z, \W) \propto f_k(\Z, \W)\Omega_{\scaleto{\mathcal{P}_Z}{4pt}}(\W)$. Such penalty notably prevents a trivial solution, as $\bm{0}$, like any constant vector, is a mode of $f_k(\cdot, \W)$ for all $\W$. Also note that the prior for $\W_{\scaleto{X}{4pt}}$ only conditions attraction while the prior for $\W_{\scaleto{Z}{4pt}}$ only affects repulsion. In the present work we focus solely on deciphering the probabilistic model that accounts for neighbor embedding loss functions and refer to \cite{bohm2020unifying} for a quantitative study of attraction and repulsion in these methods.

\paragraph{Global Structure Preservation. }
To gain intuition, consider that $\W_{\scaleto{X}{4pt}}$ is observed. As we showed in \cref{sec:within_CC}, when one relies on shift invariant kernels, the positions of the CC means are taken from a diffuse distribution. Since the above methods are all derived from the limit posteriors at $\varepsilon \to 0$, $\X_{\scaleto{M}{4pt}}$ and $\Z_{\scaleto{M}{4pt}}$ have no influence on the coupling objective. Hence if two nodes belong to different CCs, their low dimensional pairwise distance will likely not be faithful. We can expect this phenomenon to persist when the expectation on $\W_{\scaleto{X}{4pt}}$ is considered, especially when clusters are well distinguishable in $\X$. This observation is central to understand the large scale deficiency of these methods. Note that this happens at the benefit of the local structure which is faithfully represented in low dimension, as discussed in \cref{intro}. In the following section we propose to mitigate the global structure deficiency with non-degenerate MRF models.

\section{Towards Capturing Large-Scale Dependencies}\label{sec:towards_large_scale}

In this section, we investigate the ability of graph coupling to faithfully represent global structure in low dimension. To gain intuition on the case where the distribution induced by the graph is not degenerate, we consider a proper Gaussian graph coupling model and show its equivalence with PCA. We then provide a new initialization procedure to alleviate the large-scale deficiency of graph coupling when degenerate MRFs are used.

\subsection{PCA as Graph Coupling}

As we argue that the inability of SNE-like methods to reproduce the coarse-grain dependencies of the input in the latent space is due to the degeneracy of the conditional (\ref{eq:proba_perp}), a natural solution would be to consider graphical models that are well defined and integrable on the entire definition spaces of $\X$ and $\Z$. For simplicity, we consider the Gaussian model and leave the extension to other kernels for future works. Note that in this case, integrability translates into the precision matrix being full-rank. As we see with the following, the natural extension of our framework to such models leads to a well-established PCA algorithm. In the following, for a continuous variable $\bm{\Theta}_{\scaleto{Z}{4pt}}$, $\mathbb{P}(\bm{\Theta}_{\scaleto{Z}{4pt}} = \cdot)$ denotes its density.
\begin{theorem}\label{PCA_graph_coupling}
Let $\nu \geq n$,  $\bm{\Theta}_{\scaleto{X}{4pt}} \sim \mathcal{W}(\nu, \bm{I}_n)$ and $\bm{\Theta}_{\scaleto{Z}{4pt}} \sim \mathcal{W}(\nu + p - q, \bm{I}_n)$. Assume that $\bm{\Theta}_{\scaleto{X}{4pt}}$ and $\bm{\Theta}_{\scaleto{Z}{4pt}}$ structure the rows of respectively $\X$ and $\Z$ such that: 
\begin{align}
    \mathrm{vec}(\X) | \bm{\Theta}_{\scaleto{X}{4pt}} &\sim \mathcal{N}(\bm{0}, \bm{\Theta}_{\scaleto{X}{4pt}}^{-1} \otimes \bm{I}_p), \label{eq:X_given_theta} \\
    \mathrm{vec}(\Z) | \bm{\Theta}_{\scaleto{Z}{4pt}} &\sim \mathcal{N}(\bm{0}, \bm{\Theta}_{\scaleto{Z}{4pt}}^{-1} \otimes \bm{I}_q) \label{eq:Z_given_theta} \:.
\end{align}
Then the solution of the precision coupling problem:
\begin{align*}
    \min_{\Z \in \mathbb{R}^{n \times q}} -\mathbb{E}_{\bm{\Theta}_{\scaleto{X}{3pt}} | \X}\left[\log \mathbb{P}(\bm{\Theta}_{\scaleto{Z}{4pt}}=\bm{\Theta}_{\scaleto{X}{4pt}}|\Z)\right]
\end{align*}
is a PCA embedding of $\X$ with $q$ components.
\end{theorem}
We now highlight the parallels with the previous construction done for neighbor embedding methods. First note that the multivariate Gaussian with full-rank precision is inherently a pairwise MRF \cite{rue2005gaussian}. When choosing the Gaussian kernel for neighbor embedding methods, we saw that the graph Laplacian $\bm{L}_{\scaleto{X}{4pt}}$ of $\W_{\scaleto{X}{4pt}}$ was playing the role of the among-row precision matrix, as we had $\X | \W_{\scaleto{X}{4pt}} \sim \mathcal{N}(\bm{0}, \bm{L}_{\scaleto{X}{4pt}}^{-1} \otimes \bm{I}_p)$ (equation \ref{eq:gaussian_kernel}). Recall that the latter always has a null space which is spanned by the CC indicator vectors of $\W$ (\cref{sec:laplacian_prop}). Here, the key difference is that we impose a full-rank constraint on the precision $\bm{\Theta}$. Concerning the priors, we choose the ones that are conjugate to the conditionals (\ref{eq:X_given_theta}) and (\ref{eq:Z_given_theta}), as previously done when constructing the prior for neighbor embedding methods (definition \ref{def:prior_W}). Hence in the full-rank setting, the prior simply amounts to a Wishart distribution denoted by $\mathcal{W}$.

The above theorem further highlights the flexibility and generality of the graph coupling framework. Unlike usual constructions of PCA or probabilistic PCA \cite{tipping1999probabilistic}, in the above the linear relation between $\X$ and $\Z$ is recovered by solving the graph coupling problem and not explicitly stated beforehand. To the best of our knowledge, it is the first time such a link is uncovered between PCA and SNE-like methods. In contrast with the latter, PCA is well-known for its ability to preserve global structure while being significantly less efficient at identifying clusters \cite{anowar2021conceptual}. Therefore, as suspected in \cref{sec:interpretations}, the degeneracy of the conditional distribution given the graph is key to determining the distance preservation properties of the embeddings. We propose in \cref{sec:hierarchical_modelling} to combine both graph coupling approaches to strike a balance between global and local structure preservation.

\subsection{Hierarchical Graph Coupling}\label{sec:hierarchical_modelling}

The goal of this section is to show that global structure in SNE-like embeddings can be improved by structuring the CCs' positions. We consider the following hierarchical model for $\X$, where $\mathcal{P}_{\scaleto{X}{4pt}} \in \{B,D,E\}$, $k_x$ satisfies the assumptions of \cref{prop:integrability_pairwise_MRF} and $\nu_{\scaleto{X}{4pt}} \geq n$:
\begin{align*}
    \W_{\scaleto{X}{4pt}} \sim \mathbb{P}_{\scaleto{\mathcal{P}_{\scaleto{X}{3pt}},k_x}{5pt}}^{\varepsilon}(\cdot \: ; \bm{1},1), &\quad \bm{\Theta}_{\scaleto{X}{4pt}} | \W_{\scaleto{X}{4pt}} \sim \mathcal{W}(\nu_{\scaleto{X}{4pt}}, \bm{I}_{\scaleto{R}{4pt}}) \\
    \X_{\scaleto{C}{4pt}} |\W_{\scaleto{X}{4pt}} \sim \mathbb{P}_{k_x}(\cdot \:| \W_{\scaleto{X}{4pt}}), &\quad \mathrm{vec}(\X_{\scaleto{M}{4pt}}) | \bm{\Theta}_{\scaleto{X}{4pt}} \sim \mathcal{N}\left(\bm{0}, \left(\varepsilon \bm{U}_{\scaleto{[:R]}{5pt}}  \bm{\Theta}_{\scaleto{X}{4pt}}\bm{U}^\top_{\scaleto{[R]}{5pt}}\right)^{-1} \otimes \bm{I}_p\right)
\end{align*}
where $\bm{U}_{\scaleto{[R]}{5pt}}$ are the eigenvectors associated to the Laplacian null-space of $\overline{\W}_{\scaleto{X}{4pt}}$. Given a graph $\W_{\scaleto{X}{4pt}}$, the idea is to structure the CCs' relative positions with a full-rank Gaussian model.
The same model is considered for $\W_{\scaleto{Z}{4pt}}$, $\bm{\Theta}_{\scaleto{Z}{4pt}}$ and $\Z$, choosing $\nu_{\scaleto{Z}{4pt}} = \nu_{\scaleto{X}{4pt}} + p - q$ for the Wishart prior to satisfy the assumption of \cref{PCA_graph_coupling}.  With this in place, we aim at providing a complete coupling objective, matching the pairs  $(\W_{\scaleto{X}{4pt}},\bm{\Theta}_{\scaleto{X}{4pt}})$ and  $(\W_{\scaleto{Z}{4pt}},\bm{\Theta}_{\scaleto{Z}{4pt}})$. The joint negative cross-entropy can be decomposed as follows:
\begin{align}
    &\mathbb{E}_{(\W_{\scaleto{X}{3pt}}, \bm{\Theta}_{\scaleto{X}{3pt}})|\X}\left[\log \mathbb{P}((\W_{\scaleto{Z}{4pt}},\bm{\Theta}_{\scaleto{Z}{4pt}}) = (\W_{\scaleto{X}{4pt}},\bm{\Theta}_{\scaleto{X}{4pt}}) | \Z)\right] \nonumber\\
    &= \mathbb{E}_{\W_{\scaleto{X}{3pt}}|\X}\left[\log \mathbb{P}(\W_{\scaleto{Z}{4pt}} = \W_{\scaleto{X}{4pt}} | \Z)\right] + \label{eq:loss_LW} \\
    & \mathbb{E}_{(\W_{\scaleto{X}{3pt}},\bm{\Theta}_{\scaleto{X}{3pt}})|\X}\left[ \log \mathbb{P}(\bm{\Theta}_{\scaleto{Z}{4pt}} = \bm{\Theta}_{\scaleto{X}{4pt}}| \W_{\scaleto{Z}{4pt}} = \W_{\scaleto{X}{4pt}}, \Z) \right] \label{eq:add_term_Coupling}
\end{align}
where (\ref{eq:loss_LW}) is the usual coupling criterion of $\W_X$ and $\W_Z$ capturing intra-CC variability while (\ref{eq:add_term_Coupling}) is a penalty resulting from the Gaussian structure on $\mathcal{S}_{\scaleto{M}{4pt}}$. Constructed as such, the above objective allows a trade-of between local and global structure preservation. Following current trends in DR \cite{kobak2021initialization}, we propose to take care of the global structure first \textit{i.e.}\ focusing on (\ref{eq:add_term_Coupling}) before (\ref{eq:loss_LW}). The difficulty of dealing with (\ref{eq:add_term_Coupling}) lies in the hierarchical construction of the graph and the Gaussian precision (see \cref{fig:graphical_model_hierarchical}). We state the following result.

\begin{wrapfigure}[15]{R}{0.5\textwidth}
\begin{center}
\centerline{\includegraphics[width=0.5\columnwidth]{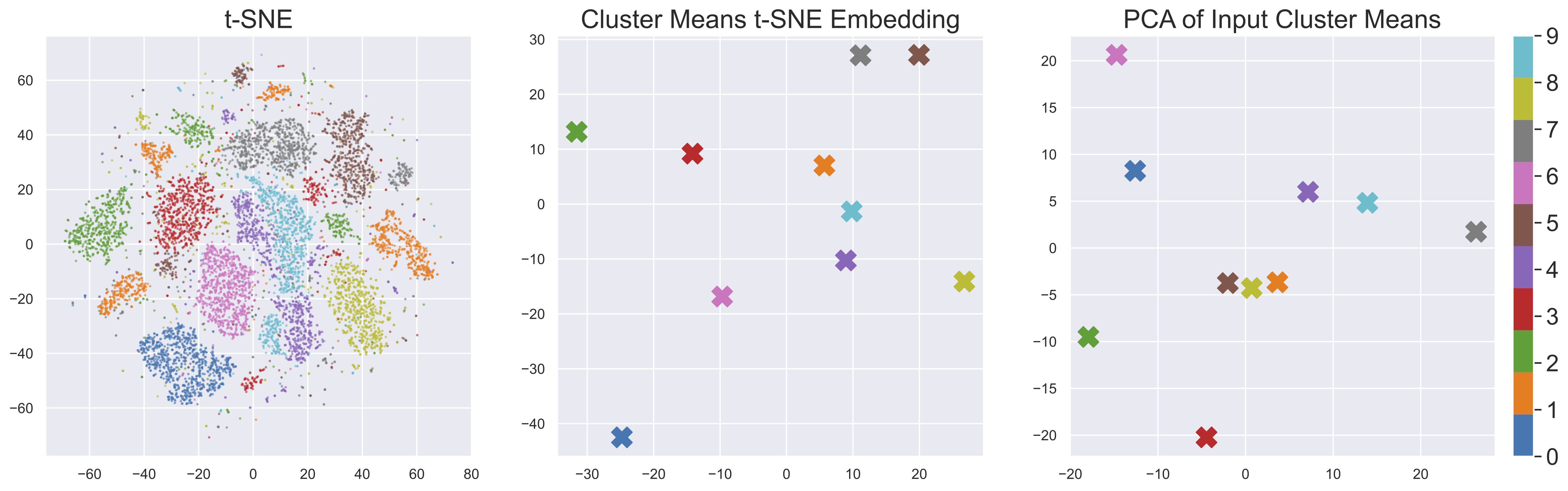}}
\caption{Left: MNIST t-SNE (perp : 30) embeddings initialized with i.i.d $\mathcal{N}(0,1)$ coordinates. Middle: using these t-SNE embeddings, mean coordinates for each digit are represented. Right: we compute a matrix of mean input coordinates for each of the $10$ digits and embed it using PCA. For t-SNE embeddings, the positions of clusters vary across different runs and don't visually match the PCA embeddings of input mean vectors (right plot).}
\label{fig:tSNE-clusters-truth}
\end{center}
\end{wrapfigure}

\begin{corollary}\label{corollary_ccPCA}
Let $\W_{\scaleto{X}{4pt}} \in \mathcal{S}_{\scaleto{W}{4pt}}$, $\bm{L} = L(\overline{\W}_{\scaleto{X}{4pt}})$ and $\mathcal{S}^q_{\scaleto{M}{4pt}}= (\ker \bm{L}) \otimes \mathbb{R}^q$, then for all $\varepsilon > 0$, given the above hierarchical model, the solution of the problem:
$$\min_{\Z \in \mathcal{S}^q_{M}} \: -\mathbb{E}_{\bm{\Theta}_{\scaleto{X}{3pt}}| \X}\left[ \log \mathbb{P}(\bm{\Theta}_{\scaleto{Z}{4pt}} = \bm{\Theta}_{\scaleto{X}{4pt}}| \W_{\scaleto{Z}{4pt}} = \W_{\scaleto{X}{4pt}}, \Z) \right]$$
is a PCA embedding of $\bm{U}_{\scaleto{[:R]}{5pt}}\bm{U}_{\scaleto{[R]}{5pt}}^\top\X$ where $\bm{U}_{\scaleto{[:R]}{5pt}}$ are the CCs' membership vectors of $\overline{\W}_{\scaleto{X}{4pt}}$.
\end{corollary}

\begin{remark}
Note that while (\ref{eq:loss_LW}) approximates the objective of SNE-like methods when $\varepsilon \to 0$, the minimizer of (\ref{eq:add_term_Coupling}) given by \cref{corollary_ccPCA} is stable for all $\varepsilon$.
\end{remark}

From this observation, we propose a simple heuristic to minimize (\ref{eq:add_term_Coupling}) that consists in computing a PCA embedding of $\mathbb{E}_{\mathbb{P}_{\scaleto{\mathcal{P}_X}{4pt}}(\cdot;\bm{K}_{\scaleto{X}{3pt}})}\left[ \bm{U}_{\scaleto{[:R]}{5pt}}\bm{U}_{\scaleto{[R]}{5pt}}^\top \right]\X$. The distribution of the connected components of the posterior of $\W_{\scaleto{X}{4pt}}$ being intractable, we resort to a Monte-Carlo estimation of the above expectation. The latter procedure called \textit{ccPCA} aims at recovering the inter-CC structure that is filtered by SNE-like methods. \textit{ccPCA} may then be used as initialization for optimizing (\ref{eq:loss_LW}) which is done by running the DR method corresponding to the graph priors at hand (\cref{sec:retrieving_DR_methods}). This second step essentially consists of refining the intra-CC structure. 

\subsection{Experiments with \textit{ccPCA}}\label{sec:ccPCA}

\begin{figure*}[t]
\begin{center}
\centerline{\includegraphics[width=\columnwidth]{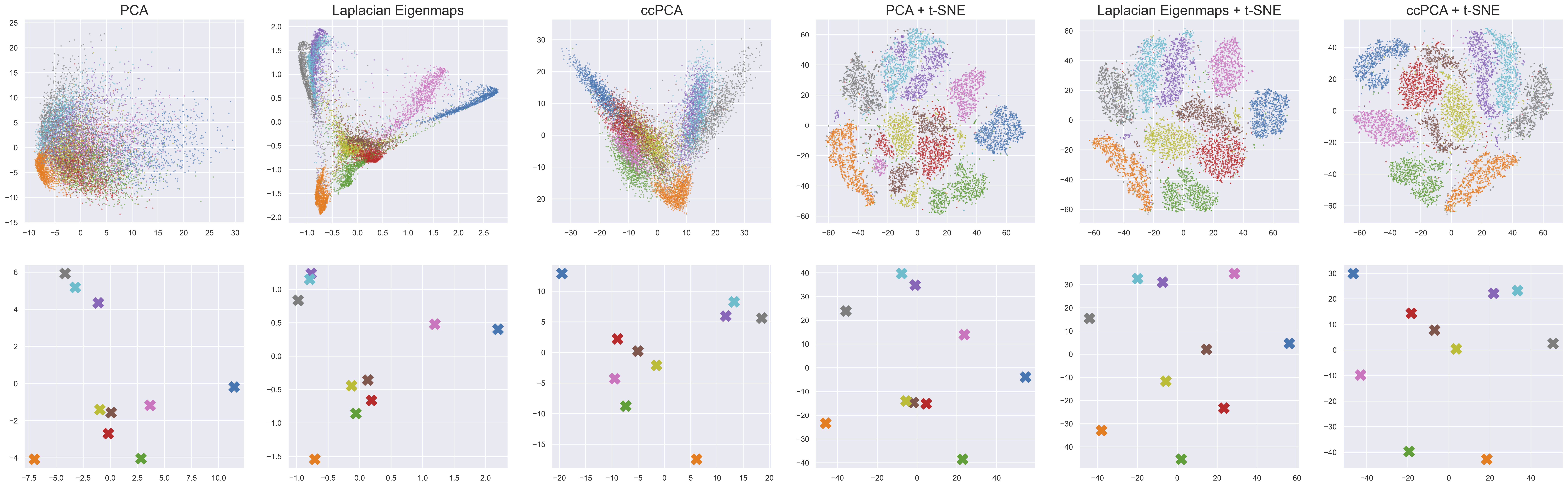}}
\caption{Top: MNIST embeddings produced by PCA, Laplacian eigenmaps, \textit{ccPCA} and finally t-SNE launched after the previous three embeddings to improve the fine-grain structure. Bottom: mean coordinates for each digit using the embeddings of the first row. The color legend is the same as in \cref{fig:tSNE-clusters-truth}. t-SNE was trained during $1000$ iterations using default parameters with the openTSNE implementation \cite{polivcar2019opentsne}.}
\label{fig:methods_embeddings}
\end{center}
\vspace{-0.8cm}
\end{figure*}

\Cref{fig:tSNE-clusters-truth} shows that a t-SNE embedding of a balanced MNIST dataset of 10000 samples \cite{deng2012mnist} with isotropic Gaussian initialization performs poorly in conserving the relative positions of clusters. As each digit cluster contains approximately $1000$ points, with a perplexity of $30$, sampling an edge across digit clusters in the graph posterior $\mathbb{P}_{\scaleto{\mathcal{P}_X}{4pt}}(\cdot;\bm{K}_{\scaleto{X}{4pt}})$ is very unlikely. Recall that the perplexity value \cite{maaten2008tSNE} corresponds to the approximate number of effective neighbors of each point. Hence images of different digits are with very high probability in different CCs of the graph posterior and their CC-wise means are not coupled as discussed in \cref{sec:interpretations}. To remedy this in practice, PCA or Laplacian eigenmaps are usually used as initialization \cite{kobak2021initialization}. 

These strategies are tested (\cref{fig:methods_embeddings}) together with \textit{ccPCA}. This shows that 
\textit{ccPCA} manages to retrieve the digits that mostly support the large-scale variability as measured by the peripheral positioning of digits $0$ (blue), $2$ (green), $6$ (pink) and $7$ (grey) given by the right side of \cref{fig:tSNE-clusters-truth}. Other perplexity values for \textit{ccPCA} are explored in appendix \ref{sec:other_perp} while the experimental setup is detailed in appendix \ref{sec:setup_exp}. In appendix \ref{sec:quantitative_evaluation}, we perform quantitative evaluations of \textit{ccPCA} for both t-SNE and UMAP on various datasets using K-ary neighborhood criteria. We find that using \textit{ccPCA} as initialization is in general more reliable than PCA and Laplacian eigenmaps for preserving global structure using both t-SNE and UMAP. 

Compared to PCA, \textit{ccPCA} manages to aggregate points into clusters, thus filtering the intra-cluster variability and focusing solely on the inter-cluster structure. Compared to Laplacian eigenmaps which perform well at identifying clusters but suffer from the same deficiency as t-SNE for positioning them, \textit{ccPCA} retains more of the coarse-grain structure. These observations support our unifying probabilistic framework and the theoretical results about MRF degeneracy which are the leading contributions of this article. The \textit{ccPCA} initialization appears as a first stepping stone towards more grounded DR methods based on the probabilistic model presented in this article.

\section{Conclusion and Perspectives}\label{Perspectives}

In this work, we shed new light on the most popular DR methods by showing that they can be unified within a common probabilistic model in the form of latent Markov Random Fields Graphs coupled by a cross-entropy. The definition of such a model constitutes a major step towards the understanding of common dimension reduction methods, in particular their structure preservation properties as discussed in this article. 

Our work offers many perspectives, among which the possibility to enrich the probabilistic model with more suited graph priors. Currently considered priors are simply the ones that are conjugate to the MRFs thus they are mostly designed to yield a tractable coupling objective. However they may not be optimal and could be modified to capture targeted features, \textit{e.g.}\ communities, in the input data, and give adapted representations in the latent space. The graph coupling approach could also be extended to more general latent structures governing the joint distribution of observations.
Finally, the probabilistic model could be leveraged to tackle hyper-parameter calibration, especially kernel bandwidths that have a great influence on the quality of the representations and are currently tuned using heuristics with unclear motivations.

\paragraph{Acknowledgments.} The authors would like to thank the anonymous reviewers whose comments and questions helped improve the clarity of this manuscript, as well as Aurélien Garivier, Antoine Barrier and Floshi Poshi for helpful discussions. This work was supported by the Agence Nationale de la Recherche  ANR-18-CE45-0023 SingleStatOmics.

\newpage

\appendix
\onecolumn

\section{Proofs}\label{sec:proofs}

\subsection{Proof of \cref{prop:integrability_pairwise_MRF}} \label{proof:lambda_perp_integrability}

$\W \in \mathcal{S}_{\scaleto{W}{4pt}}$ is the weight matrix of a graph with $R$ connected components $\{C_1, ..., C_R\}$ partitioning $[n]$. Since $k$ is upper bounded by a constant, there exists $M_+ > 1$ that upper bounds $k$. Let $\bm{\mathcal{T}}$ be the adjacency matrix of a spanning forest of $\W$, since each edge of $\W$ is bounded by $n$, one has:
\begin{align}
    \int f_{k}(\X, \W) \lambda_{\mathcal{S}_{\scaleto{C}{3pt}}}(d\X) &= \int \prod_{(i,j) \in [n]^2} k(\X_{i} - \X_{j})^{W_{ij}} \lambda_{\mathcal{S}_{\scaleto{C}{3pt}}}(d\X) \nonumber \\
    &\leq M_+^{n^{3}} \int \prod_{(i,j) \in [n]^2} k(\X_{i} - \X_{j})^{\mathcal{T}_{ij}} \lambda_{\mathcal{S}_{\scaleto{C}{3pt}}}(d\X) \nonumber \\
    &\leq M_+^{n^{3}} \prod_{r \in [R]} \int \prod_{(i,j) \in C_{r}^2} k(\X_{i} - \X_{j})^{\mathcal{T}_{ij}} \lambda_{\mathcal{S}_{\scaleto{C}{3pt}}}(d\X) \:. \label{bound_MRF_product_CC}
\end{align}
Let $r \in [R]$. The spanning tree corresponding to the $r^{th}$ connected component called $\bm{\mathcal{T}}^r$ has exactly $n_r-1$ edges. There exists a leaf node $\ell \in [n]$ of $\bm{\mathcal{T}}^r$ and let $\tilde{\ell}$ be the node linked to it. Consider a bijective map $\sigma \colon C_r \backslash \{\ell\} \to [n_r - 1]$ such that $\sigma(\tilde{\ell}) = 1$ and for $(i,j) \in (C_r \backslash \{\ell\})^2$, $\sigma(i) \leq \sigma(j)$ implies that node $i$ has a shorter path on $\overline{\bm{\mathcal{T}}^r}$\footnote{Symmetrized version \textit{i.e.}\ $\overline{\bm{\mathcal{T}}^r} = \bm{\mathcal{T}}^r + \bm{(\mathcal{T}}^r)^\top$.} to $\ell$ than node $j$. There exists a bijective map $e \colon [2:n_r - 1] \to [n_r - 2]$ such that for $i \in [2:n_R-1]$, $\overline{\bm{\mathcal{T}}^r}_{\sigma^{-1}(i), \sigma^{-1}(e(i))} > 0$ and node $\sigma^{-1}(e(i))$ has a shorter path on $\overline{\bm{\mathcal{T}}^r}$ to node $\ell$ than node $\sigma^{-1}(i)$.

Recall that since $\X \in \mathcal{S}_{\scaleto{C}{4pt}}$ one has: $\sum_{i \in C_r} \X_i = 0$ hence $\X_{\ell} = - \sum_{i \neq \ell} \X_i$. Let us now consider the linear map $\phi^r$ such that:
\begin{align*}
\forall i \in [n_r - 1], \quad \phi^r(\bm{X}_i) = \left\{
    \begin{array}{ll}
        \X_{\sigma^{-1}(i)} + \sum_{j \in [n_r - 1]} \X_{\sigma^{-1}(j)} & \mbox{if i = 1}\\
        \X_{\sigma^{-1}(i)} - \X_{\sigma^{-1}(e(i))} & \mbox{otherwise} \:.
    \end{array}
\right.
\end{align*}

We now show that the change of variable $\phi^r$ is a $\mathcal{C}^1$ diffeomorphism by proving that its Jacobian has full rank. Ordering the columns with the map $\sigma$, the latter takes the form:
\[
    \bm{J}_{\phi^r} = \left(
    \begin{array}{ccccc}
    2 & 1 & 1 & \dots & 1 \\
      & 1 & 0 & \dots & 0 \\
      &   & \ddots &  \ddots & \vdots \\
      & \bm{A} &   & \ddots & 0 \\
      &               &   &   & 1
    \end{array}
    \right)
\]
where $\bm{A}$ is a strictly lower triangular matrix such that for all $i \in [2:n_r-1]$, $A_{ie(i)} = -1$ and for all $t \neq e(i)$, $A_{it}=0$. The above can be factorized as:
\[
\bm{J}_{\phi^r} = 
\left(
    \begin{array}{ccccc}
    \alpha_{n_r-1} & \alpha_{n_r-2} & \dots & \alpha_2 & \alpha_1 \\
    0  & 1 & 0 & \dots & 0 \\
    \vdots & \ddots & \ddots & \ddots & \vdots \\
    \vdots & & \ddots & \ddots & 0 \\
    0 & \dots & \dots & 0 & 1
    \end{array}
    \right)^{-1}
\left(
    \begin{array}{ccccc}
    1 & 0 & \dots & \dots & 0 \\
      & 1 & \ddots & & \vdots\\
      & & \ddots & \ddots & \vdots \\
      & \bm{A} & & \ddots & 0 \\
      & & & & 1
    \end{array}
\right)
\]
where $\alpha_{1}=-1$ and for $\ell > 1$, $\alpha_\ell = \sum_{j < l} \alpha_j\ind_{e(n_r - j)=n_r -\ell} - 1$. With this in place, for $i \in [n_r -1]$, $\alpha_i \neq 0$ in particular $\alpha_{n_r-1} \neq 0$ therefore $|\bm{J}_{\phi^r}| \neq 0 $ and $\phi^r$ is a $\mathcal{C}^1$ diffeomorphism. This change of variable yields:
\begin{align*}
\int \prod_{(i,j) \in C_{r}^2} k(\X_{i} - \X_{j})^{\mathcal{T}_{ij}} \lambda_{\mathcal{S}_{\scaleto{C}{3pt}}}(d\X) 
&= \int \bigotimes_{i \in [n_r - 1]} k(\bm{Y}_i) |\bm{J}_{\phi^r}(\bm{Y})|^{-1} \lambda_{\mathbb{R}^p}(d\bm{Y}) \\
&= |\bm{J}_{\phi^r}|^{-1} \prod_{i \in [n_r - 1]} \int k(\bm{Y}_i) \lambda_{\mathbb{R}^p}(d\bm{Y}_i)
\end{align*}
using the Fubini Tonelli theorem. The result follows from $\lambda_{\mathbb{R}^p}$-integrability of $k$ and upper bound \ref{bound_MRF_product_CC}.

\subsection{Proof of \cref{prop:posterior_W}}
\label{proof:posterior_limit}

Let $\mathcal{P} \in \{B, D, E\}$, $k$ be a valid kernel (assumptions of \cref{prop:integrability_pairwise_MRF}) with $\bm{K}_{\scaleto{X}{4pt}} = (k(\X_{i} - \X_{j}))_{(i,j) \in [n]^2}$ and $\bm{\pi} \in \mathbb{R}_+^{n \times n}$. Let $\W \sim \mathbb{P}_{\scaleto{\mathcal{P},k}{5pt}}^{\varepsilon}(\cdot \: ; \bm{\pi},1)$. Inversion of conditional with Bayes rule gives:
\begin{align}
    \forall \bm{W} \in \mathcal{S}_{\scaleto{W}{4pt}}, \quad \mathbb{P}(\W | \bm{X}) \propto
    \mathcal{C}_{k}^{\varepsilon}(\W)^{-1} f^{\varepsilon}(\X, \W) f_{k}(\X, \W) \mathbb{P}^{\varepsilon}_{\scaleto{\mathcal{P},k}{5pt}}(\W; \bm{\pi}, 1) \label{inversion_Conditional}
\end{align}
where the prior reads:
\begin{align}
    \mathbb{P}_{\scaleto{\mathcal{P},k}{5pt}}^{\varepsilon}(\bm{W}; \bm{\pi}, 1) \propto \mathcal{C}^{\varepsilon}_k(\W) \Omega_{\scaleto{\mathcal{P}}{4pt}}(\W) \prod_{(i,j) \in [n]^2} \pi_{ij}^{W_{ij}} \:.
\end{align}
Hence the joint normalizing constant simplifies such that:
\begin{align}
    \forall \bm{W} \in \mathcal{S}_{\scaleto{W}{4pt}}, \quad \mathbb{P}(\W | \bm{X}) &\propto
    f^{\varepsilon}(\X, \W) \Omega_{\scaleto{\mathcal{P}}{4pt}}(\W) \prod_{(i,j) \in [n]^2} \left(\pi_{ij} k(\X_i - \X_j)\right)^{W_{ij}} \\
    &\xrightarrow[\varepsilon \to 0]{} \Omega_{\scaleto{\mathcal{P}}{4pt}}(\W) \prod_{(i,j) \in [n]^2}  \left(\pi_{ij} k(\X_i - \X_j)\right)^{W_{ij}}
\end{align}
which ends the proof. As a complement, we now explicit the simple forms taken by the posterior limit graph in each case.

\paragraph{$B$-Prior}
Recall that in this case the prior reads:
\begin{align*}
    \mathbb{P}_{\scaleto{B}{4pt}}^{\varepsilon}(\bm{W}; \bm{\pi},1) &\propto \mathcal{C}_{k}^{\varepsilon}(\bm{W}) \prod_{(i,j) \in [n]^2} \pi_{ij}^{W_{ij}} \ind_{W_{ij} \leq 1} \:.
\end{align*}
Therefore the posterior limit graph has the distribution:
\begin{align*}
    \mathbb{P}_{\scaleto{B}{4pt}}(\W ;\bm{\pi} \odot \bm{K}_{\scaleto{X}{4pt}})
    &= \frac{\prod_{(i,j) \in [n]^2}  \left(\pi_{ij} k(\X_i - \X_j)\right)^{W_{ij}} \ind_{W_{ij} \leq 1}}{\sum_{\W \in \mathcal{S}_{\scaleto{W}{4pt}}} \prod_{(i,j) \in [n]^2}  \left(\pi_{ij} k(\X_i - \X_j)\right)^{W_{ij}} \ind_{W_{ij} \leq 1}} \\
    &= \prod_{(i,j) \in [n]^2}  \left(\frac{\pi_{ij} k(\X_i - \X_j)}{1+\pi_{ij} k(\X_i - \X_j)}\right)^{W_{ij}} \left(\frac{1}{1+\pi_{ij} k(\X_i - \X_j)}\right)^{1-W_{ij}} \ind_{W_{ij} \leq 1} \:.
\end{align*}

This distribution amounts to: $\forall (i,j) \in [n]^2, \quad \W_{ij} \stackrel{\perp\!\!\!\!\perp}{\sim} \mathcal{B}\left( \frac{\pi_{ij} k(\X_i - \X_j)}{1+\pi_{ij} k(\X_i - \X_j)} \right)$.

\paragraph{$D$-Prior} The prior writes:
\begin{align*}
    \mathbb{P}_{\scaleto{D}{4pt}}^{\varepsilon}(\bm{W}; \bm{\pi}, 1) &\propto \mathcal{C}_{k}^{\varepsilon}(\bm{W}) \prod_{(i,j) \in [n]^2} \pi_{ij}^{W_{ij}} \ind_{W_{i+} = 1} \:.
\end{align*}
The distribution of the posterior limit then becomes:
\begin{align*}
    \mathbb{P}_{\scaleto{D}{4pt}}(\W ;\bm{\pi} \odot \bm{K}_{\scaleto{X}{4pt}}) &= \frac{\prod_{(i,j) \in [n]^2}  \left(\pi_{ij} k(\X_i - \X_j)\right)^{W_{ij}} \ind_{W_{i+} = 1}}{\sum_{\W \in \mathcal{S}_{\scaleto{W}{4pt}}} \prod_{(i,j) \in [n]^2}  \left(\pi_{ij} k(\X_i - \X_j)\right)^{W_{ij}} \ind_{W_{i+} = 1}} \\
    &= \frac{\prod_{(i,j) \in [n]^2}  \left(\pi_{ij} k(\X_i - \X_j)\right)^{W_{ij}} \ind_{W_{i+} = 1}}{\prod_{i \in [n]} \sum_{\ell \in [n]} \pi_{i\ell} k(\X_i - \X_\ell)} \\
    &= \prod_{(i,j) \in [n]^2} \left(\frac{\pi_{ij} k(\X_i - \X_j)}{\sum_{\ell \in [n]} \pi_{i\ell} k(\X_i - \X_\ell)}\right)^{W_{ij}} \ind_{W_{i+} = 1} \:.
\end{align*}

This distribution amounts to: $\forall i \in [n], \quad \W_{i} \stackrel{\perp\!\!\!\!\perp}{\sim} \mathcal{M}\left(1, \left(\frac{\pi_{ij} k(\X_i - \X_j)}{\sum_{\ell \in [n]} \pi_{i\ell} k(\X_i - \X_\ell)}\right)_{j \in [n]}\right)$.

\paragraph{$E$-Prior}
In this case the prior reads:
\begin{align*}
    \mathbb{P}_{\scaleto{E}{4pt}}^{\varepsilon}(\bm{W}; \bm{\pi}, 1) &\propto \mathcal{C}_{k}^{\varepsilon}(\bm{W}) \prod_{(i,j) \in [n]^2} \frac{\pi_{ij}^{W_{ij}}}{W_{ij}!} \ind_{W_{++} = n} \:.
\end{align*}
Finally, deriving the distribution of the posterior graph limit:
\begin{align*}
    \mathbb{P}_{\scaleto{E}{4pt}}(\W ;\bm{\pi} \odot \bm{K}_{\scaleto{X}{4pt}}) &= \frac{\prod_{(i,j) \in [n]^2}  (W_{ij}!)^{-1}\left(\pi_{ij} k(\X_i - \X_j)\right)^{W_{ij}} \ind_{W_{++} = n}}{\sum_{\W \in \mathcal{S}_{\scaleto{W}{4pt}}} \prod_{(i,j) \in [n]^2} (W_{ij}!)^{-1} \left(\pi_{ij} k(\X_i - \X_j)\right)^{W_{ij}} \ind_{W_{++} = n}} \\
    &= n! \prod_{(i,j) \in [n]^2} (W_{ij})^{-1} \left(\frac{\pi_{ij} k(\X_i - \X_j)}{\sum_{(\ell,t) \in [n]^2} \pi_{\ell t} k(\X_\ell - \X_t)}\right)^{W_{ij}} \ind_{W_{++} = n} \:.
\end{align*}

This distribution amounts to: $\W \sim \mathcal{M}\left(n, \left(\frac{\pi_{ij} k(\X_i - \X_j)}{\sum_{(\ell,t) \in [n]^2} \pi_{\ell t} k(\X_\ell - \X_t)}\right)_{(i,j) \in [n]^2}\right)$.

\subsection{Proof of \cref{PCA_graph_coupling}}

We consider the following hierarchical model, for $\nu_{\scaleto{X}{4pt}}, \nu_{\scaleto{Z}{4pt}} \geq n$:
\begin{align*}
    \bm{\Theta}_{\scaleto{X}{4pt}} &\sim  \mathcal{W}(\nu_{\scaleto{X}{4pt}}, \bm{I}_n) \\
    \mathrm{vec}(\X) | \bm{\Theta}_{\scaleto{X}{4pt}} &\sim \mathcal{N}(\bm{0}, \bm{\Theta}_{\scaleto{X}{4pt}}^{-1} \otimes \bm{I}_p) \\
    \bm{\Theta}_{\scaleto{Z}{4pt}} &\sim  \mathcal{W}(\nu_{\scaleto{Z}{4pt}}, \bm{I}_n) \\
    \mathrm{vec}(\Z) | \bm{\Theta}_{\scaleto{Z}{4pt}} &\sim \mathcal{N}(\bm{0}, \bm{\Theta}_{\scaleto{Z}{4pt}}^{-1} \otimes \bm{I}_q) \:.
\end{align*}
With this at hand, the posteriors for $\bm{\Theta}_X$ and $\bm{\Theta}_Z$ can be derived in closed form: 
\begin{align*}
    \bm{\Theta}_{\scaleto{X}{4pt}} | \X &\sim  \mathcal{W}(\nu_{\scaleto{X}{4pt}}+p, \left(\bm{I}_n + \X \X^\top\right)^{-1}) \\
    \bm{\Theta}_{\scaleto{Z}{4pt}} | \Z &\sim  \mathcal{W}(\nu_{\scaleto{Z}{4pt}} + q, \left(\bm{I}_n + \Z \Z^\top\right)^{-1}) \:.
\end{align*}

Keeping terms of $-\mathbb{E}_{\bm{\Theta}_{\scaleto{X}{4pt}}}[\log \mathbb{P}(\bm{\Theta}_{\scaleto{Z}{4pt}} = \bm{\Theta}_{\scaleto{X}{4pt}}| \Z )| \X]$ that depends on $\Z$, one has the optimization problem:
\begin{align*}
    \min_{\Z \in \mathbb{R}^{n \times q}} \quad \frac{\nu_{\scaleto{X}{4pt}}+p}{2}\operatorname{tr}\left(\Z^\top(\bm{I}_n +  \X\X^\top)^{-1}\Z\right) - \frac{\nu_{\scaleto{Z}{4pt}}+q}{2}\log |\bm{I}_n +  \Z\Z^\top|
\end{align*}
Our strategy is to first find the optimal sample covariance matrix $\Z\Z^\top$ and then focus on the solution in $\Z$. To that extend, consider the eigendecomposition of the sample covariance matrices: $\X\X^\top = \bm{V D} \bm{V}^\top$ and $\Z\Z^\top = \bm{U \Lambda} \bm{U}^\top$ where $\bm{D}=\operatorname{diag}(\bm{d})$ and $\bm{\Lambda}=\operatorname{diag}(\bm{\lambda})$ such that $d_1 \geq ... \geq d_n$ and $\lambda_1 \geq ... \geq \lambda_n$. Denoting $\gamma = (\nu_{\scaleto{X}{4pt}}+q)/(\nu_{\scaleto{Z}{4pt}}+p)$, we consider the following problem:
\begin{align}
   \min_{\bm{U} \in \mathcal{O}(n), \bm{\Lambda}} \quad & \operatorname{tr}\left(\bm{U} \bm{\Lambda} \bm{U}^\top \bm{V} (\bm{I}_n + \bm{D})^{-1} \bm{V}^\top\right) - \gamma \log |\bm{I}_n + \bm{\Lambda}| \label{eq:optim_eigenvalues_eigenvectors} \\
    \textrm{s.t.} \quad & \bm{\Lambda} \succcurlyeq \bm{0} \label{eq:positive_definite_constraint}\\
    & \operatorname{rank}(\bm{\Lambda}) \leq q \label{eq:rank_constraint}
\end{align}

The above problem is non-convex because of the rank constraint (\ref{eq:rank_constraint}). Nonetheless it can be simplified as we now show. 

We focus on finding the optimal eigenvectors first. To that extent, let us denote, $\bm{R} = \bm{U}^\top\bm{V}$. Only the left term in (\ref{eq:optim_eigenvalues_eigenvectors}) depends on $\bm{R}$. The optimization problem for eigenvectors writes:
\begin{align}
   \min_{\bm{R} \in \mathcal{O}(n)} \quad & \operatorname{tr}\left(\bm{R}^\top \bm{\Lambda} \bm{R} (\bm{I}_n + \bm{D})^{-1} \right) \label{eq:optim_eigenvalues_eigenvectors}
\end{align}
The objective (\ref{eq:optim_eigenvalues_eigenvectors}) can be expressed as: $\sum_{(i,j) \in [n]^2} \lambda_i (1 + d_j)^{-1} R_{ij}^2$. Now one can notice that since $\bm{R}$ is orthogonal, $\bm{R} \odot \bm{R}$ is doubly stochastic (\textit{i.e.}\ sum of coefficients on each row and column is equal to one). Therefore thanks to the Birkhoff–von Neumann theorem, there exists $\theta_1, ..., \theta_L \geq 0$, $\sum_{\ell \in [L]} \theta_\ell = 1$ and permutation matrices $\bm{P}_1, ..., \bm{P}_L$ such that:
$$\bm{R} \odot \bm{R} = \sum_{\ell \in [L]} \theta_\ell \bm{P}_\ell$$
where for all $\ell \in [L]$, there exists a permutation $\sigma_\ell$ of $[n]$ such that $P_{\ell,ij} = \ind_{\sigma_{\ell}(i) = j}$ for $(i,j) \in [n]^2$. 

With this at hand, objective (\ref{eq:optim_eigenvalues_eigenvectors}) writes: $\sum_{\ell \in [L]} \theta_\ell \sum_{i \in [n]} \lambda_i (1 + d_{\sigma_\ell(i)})^{-1}$. There exists a permutation $\sigma^\star$ such that the quantity $\sum_{i \in [n]} \lambda_i (1 + d_{\sigma_\ell(i)})^{-1}$ is minimal. Note that the identity permutation \textit{i.e.}\ for $i \in [n]$, $ \sigma(i) = i$ is optimal in this case as the $(\lambda_i)_{i \in [n]}$ and the $(d_i)_{i \in [n]}$ are in decreasing order. Then choosing for $\ell \in [L]$, $\theta_\ell = _{\sigma_\ell = \sigma^\star}$ minimizes the latter quantity. Therefore the solution of (\ref{eq:optim_eigenvalues_eigenvectors}) $\bm{R}^{\star}$ is such that for $(i,j) \in [n]^2$, $R^\star_{ij} = \pm \ind_{\sigma^\star(i)=j}$. Thus an optimum in $\bm{U}$ of $\ref{eq:optim_eigenvalues_eigenvectors}$ is such that $\bm{U}^\star = \bm{V} \bm{R}^\star$. 

Hence $\bm{U} = \bm{V}$, in particular, is optimal. We will choose this $\bm{U}$ in what follows as the sign of the axes do not influence the characterization of the final result in $\Z$ as a PCA embedding. Such a choice gives $\Z \Z^\top = \bm{V} \bm{\Lambda} \bm{V}^\top$. 

Now it remains to find the optimal eigenvalues $(\lambda_i)_{i \in [n]}$. The rank constraint (\ref{eq:rank_constraint}) can be easily dealt with: since the eigenvalues are sorted in decreasing order, the constraint implies that for $i \geq q$, $\lambda_i=0$.  Thus the eigenvalue problem can be formulated in $\mathbb{R}^q$:
\begin{align}
    \min_{\bm{\lambda} \in \mathbb{R}^q} \quad & \bm{\lambda}^\top (\bm{1} + \bm{d})^{-1} - \gamma \bm{1}^\top \log (\bm{1} + \bm{\lambda}) \label{eq:objective_lambda}\\
    \textrm{s.t.} \quad & \forall i \in [q], \quad  \lambda_i \geq 0 , \quad \lambda_1 \geq ... \geq \lambda_q \label{eq:feasibility_lambda}
\end{align}
where (\ref{eq:feasibility_lambda}) accounts for (\ref{eq:positive_definite_constraint}). The above is convex. (\ref{eq:objective_lambda}) is minimized for $\bm{\lambda} = \gamma (\bm{1} + \bm{d}) - \bm{1}$. Taking the feasibility constraint (\ref{eq:feasibility_lambda}) into account one has a solution $\bm{\lambda}^*$ such that:
$$\forall i \in [n], \quad 
\lambda_i^* = \left\{
    \begin{array}{ll}
        \max(0, \gamma(1 + d_i) - 1) \quad &\mbox{if} \quad i \leq q \\
        0 \quad &\mbox{otherwise} \:.
    \end{array}
\right. $$

Note that this solution is not unique if there are repeated eigenvalues. Notice also that one has the freedom to choose the Wishart prior parameters such that $\gamma=1$. Doing so, the solution satisfies $\Z^\star \Z^{\star \:T} = \bm{V}_{[:,q]} \bm{D}_{[q,q]} \bm{V}^\top_{[q,:]}$. Therefore there exists $\bm{R}$ an orthogonal matrix of size $q$ such that $\Z^{\star} = \bm{V}_{[:,q]}\bm{D}_{[q,q]}^{\frac{1}{2}}\bm{R}$. The latter is the output of a PCA model of $\X$ with $q$ components, which is defined up to a rotation.

\subsection{Proof of Corollary \ref{corollary_ccPCA}}

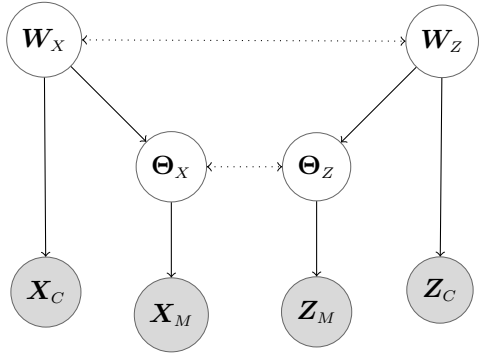
\begin{figure}
    \centering
    \begin{tikzpicture}[
    roundnode_w/.style={circle, draw=black!60, fill=white, minimum size=8.5mm, line width=0.1mm},
    roundnode_g/.style={circle, draw=black!60, fill=black!15, minimum size=8.5mm, line width=0.1mm},
    squarednode/.style={rectangle, draw=black!60, fill=black!5, minimum size=7mm}
    ]
        \node[roundnode_w]  (Theta_X)     {$\bm{\Theta}_{\scaleto{X}{4pt}}$};
    \node[roundnode_w]      (W_X)   [above left=of Theta_X]  {$\W_{\scaleto{X}{4pt}}$};
    \node[roundnode_g]        (X_C)     [below left =of Theta_X]  {$\X_{\scaleto{C}{4pt}}$};
    \node[roundnode_g]        (X_M)     [below =of Theta_X]  {$\X_{\scaleto{M}{4pt}}$};
    \node[roundnode_w]      (Theta_Z)       [right=of Theta_X]    {$\bm{\Theta}_{\scaleto{Z}{4pt}}$};
    \node[roundnode_w]      (W_Z)       [above right=of Theta_Z]    {$\W_{\scaleto{Z}{4pt}}$};
    \node[roundnode_g]        (Z_C)     [below right=of Theta_Z]  {$\Z_{\scaleto{C}{4pt}}$};
    \node[roundnode_g]        (Z_M)     [below =of Theta_Z]  {$\Z_{\scaleto{M}{4pt}}$};
    
        \draw[->] (W_X) -- (X_C);
    \draw[->] (W_X) -- (Theta_X);
    \draw[->] (Theta_X) -- (X_M);
    \draw[->] (W_Z) -- (Z_C);
    \draw[->] (W_Z) -- (Theta_Z);
    \draw[->] (Theta_Z) -- (Z_M);
    \draw[dotted,<->] (W_X) -- (W_Z);
    \draw[dotted,<->] (Theta_X) -- (Theta_Z);
    \end{tikzpicture}
    \caption{Graphical representation of the hierarchical model considered in \cref{sec:hierarchical_modelling}. Plain directed arrows represent conditional dependencies while dotted arrows represent the coupling links. \Cref{corollary_ccPCA} provides a solution for the coupling between $\bm{\Theta}_{\scaleto{X}{4pt}}$ and $\bm{\Theta}_{\scaleto{Z}{4pt}}$.}
    \label{fig:graphical_model_hierarchical}
\end{figure}

With the presented hierarchical model (\cref{fig:graphical_model_hierarchical}), the coupling problem is the following:
\begin{align}\label{eq:optim_corollary}
    \min_{\Z \in \mathcal{S}^q_{\scaleto{M}{3pt}}} \quad \operatorname{tr}\left(\bm{U}_{\scaleto{[:R]}{5pt}}\Z^\top(\bm{I}_{\scaleto{R}{4pt}} +  \varepsilon \bm{U}_{\scaleto{[R]}{5pt}}^\top\X\X^\top\bm{U}_{\scaleto{[:R]}{5pt}})^{-1}\bm{U}_{\scaleto{[R]}{5pt}}^\top\Z\right) - \log |\bm{I}_{\scaleto{R}{4pt}}  +  \varepsilon \bm{U}_{\scaleto{[R]}{5pt}}^\top\Z\Z^\top\bm{U}_{\scaleto{[:R]}{5pt}}| 
\end{align}
where $\bm{U}_{\scaleto{[:R]}{5pt}}$ are the eigenvectors associated to the Laplacian null-space of $\overline{\W}_{\scaleto{X}{4pt}}$.

Let us denote $\bar{\Z} = \bm{U}_{\scaleto{[R]}{5pt}}^\top\Z \in \mathbb{R}^{R \times q}$ and $\bar{\X} = \bm{U}_{\scaleto{[R]}{5pt}}^\top\X \in \mathbb{R}^{R \times p}$. Note that $\Z \to \bm{U}_{\scaleto{[R]}{5pt}}^\top\Z$ is a bijective linear map from $\mathcal{S}^q_{\scaleto{M}{4pt}}$ to $\mathbb{R}^{R \times q}$ with inverse $\bar{\Z} \to \bm{U}_{\scaleto{[R]}{5pt}}\bar{\Z}$ (and equivalently for $\mathbb{R}^{R \times p}$). Hence (\ref{eq:optim_corollary}) is equivalent to:
\begin{align}\label{eq:small_dim_optim_corollary}
    \min_{\bar{\Z} \in \mathbb{R}^{R \times q}} \quad \operatorname{tr}\left(\bar{\Z}^\top(\bm{I}_{\scaleto{R}{4pt}} +  \varepsilon \bar{\X}\bar{\X}^\top)^{-1}\bar{\Z}\right) - \log |\bm{I}_{\scaleto{R}{4pt}}  +  \varepsilon \bar{\Z}\bar{\Z}^\top| 
\end{align}

According to \cref{PCA_graph_coupling}, the solution of problem (\ref{eq:small_dim_optim_corollary}) is such that there exists $\bm{R}$ orthogonal, $\bar{\Z}^\star = \bm{V}_{[:,q]} \bm{S}_{[q,q]} \bm{R}$ where $\bar{\X}\bar{\X}^\top = \bm{V} \bm{S}^2 \bm{V}^\top$ is the eigendecomposition in an orthogonal basis of the among-row covariance matrix of $\bar{\X}$. Note that the solution does not depend on $\varepsilon$.

Therefore (\ref{eq:optim_corollary}) is solved for $\Z^\star = \bm{U}_{\scaleto{[:R]}{5pt}}\bm{V}_{[:,q]} \bm{S}_{[q,q]} \bm{R}$. One can notice that since the singular value decomposition (\textit{i.e.}\ SVD) of $\bm{U}_{\scaleto{[R]}{5pt}}^\top\X$ takes the form $\bm{V}\bm{S}\bm{B}$ where $\bm{B}$ is an semi-orthogonal matrix of size $p$, then $\bm{U}_{\scaleto{[:R]}{5pt}}\bm{U}_{\scaleto{[R]}{5pt}}^\top\X = \bm{U}_{\scaleto{[:R]}{5pt}}\bm{V}\bm{S}\bm{B}$. Noticing that $\bm{V}' = \bm{U}_{\scaleto{[:R]}{5pt}}\bm{V}$ is orthogonal, one has that $\bm{V}' \bm{S}\bm{B}$ is a compact SVD of $\bm{U}_{\scaleto{[:R]}{5pt}}\bm{U}_{\scaleto{[R]}{5pt}}^\top\X$. Therefore, since $\Z^\star = \bm{V}' \bm{S}$, $\Z^\star$ is a PCA embedding of $\bm{U}_{\scaleto{[:R]}{5pt}}\bm{U}_{\scaleto{[R]}{5pt}}^\top\X$.

\section{Experiments Supplementary Material}\label{sec:exp_sup_mat}

\subsection{Experimental Setup and Details About \textit{ccPCA}}\label{sec:setup_exp}

\paragraph{Implementation of existing methods.} For t-SNE, we rely on the openTSNE implementation \cite{polivcar2019opentsne} for both computing the kernel $\bm{K}_{\scaleto{X}{4pt}}$ with appropriate bandwidths and running the tSNE algorithm. We keep the training default parameters and $1000$ iterations of gradient descent. For all experiments, the default perplexity of $30$ was used to set the kernel bandwidths. For UMAP, we use the default Python implementation of \cite{mcinnes2018umap} with default parameters. For PCA and Laplacian eigenmaps, the scikit-learn implementation is used \cite{pedregosa2011scikit} with default parameters as well. 

\paragraph{\textit{ccPCA}.} The pseudo code of the algorithm is given in \cref{algo:ccPCA}. CCs' memberships (\textit{i.e.}\ eigenvectors $\bm{U}_{\scaleto{[R]}{5pt}}$) are computed using igraph \cite{csardi2006igraph}. Regarded the time complexity of ccPCA, one can sample the posterior graph with constant time if $\mathcal{P}_{\scaleto{X}{4pt}} = E$, linear time if $\mathcal{P}_{\scaleto{X}{4pt}} = D$ and quadratic time if $\mathcal{P}_{\scaleto{X}{4pt}} = B$. Moreover, computing $\bm{U}_{\scaleto{[R]}{5pt}}$ can be done with linear complexity \textit{w.r.t.}\ the number of nodes. Hence the time complexity is $O(N \times n)$ for $E$ and $D$ priors and $O(N \times n^2)$ for the $B$ prior, where $N$ is the number of Monte Carlo samples. In practice we found that $N\approx 100$ Monte Carlo samples produce a consistent \textit{ccPCA} embedding for $n \approx 10000$. Note that the time complexity of PCA is $O(\min(p^3,n^3))$ where $p$ is the dimensionality (\textit{i.e.}\ number of columns) of $\X$. Hence in most common applications involving images or biological sequencing data (where $p$ is very large), the additional time complexity brought by \textit{ccPCA} compared to PCA is negligible.

\vspace{0.5cm}

\begin{algorithm}[H]
  \caption{\textit{ccPCA}}
  \label{alg:3CPCA}
\begin{algorithmic}
  \STATE {\bfseries Input:} $\bm{K}_{\scaleto{X}{4pt}}$, $\mathcal{P}_{\scaleto{X}{4pt}}$, N
  \FOR{$\ell=1$ {\bfseries to} $N$}
  \STATE Sample $\W^{\ell} \sim  \mathbb{P}_{\scaleto{\mathcal{P}_{\scaleto{X}{4pt}}}{5pt}}(\cdot ; \bm{K}_{\scaleto{X}{4pt}})$
  \STATE Compute CCs' memberships $\bm{U}^{\ell}_{\scaleto{[R]}{5pt}}$ of $\W^{\ell}$
  \ENDFOR
  \STATE {\bfseries Output:} PCA of $\left(N^{-1}\sum_{\ell \in [N]} \bm{U}^{\ell}_{\scaleto{[R]}{5pt}} \bm{U}_{\scaleto{[R]}{5pt}}^{\ell\:T} \right) \X$
\end{algorithmic}
\label{algo:ccPCA}
\end{algorithm}

All experiments are performed on a
machine with four Intel Core i5 processors and 16 GB memory.

\subsection{\textit{ccPCA} with Varying Perplexity Values}\label{sec:other_perp}

\begin{figure*}[h]
\vskip 0.2in
\begin{center}
\centerline{\includegraphics[width=\columnwidth]{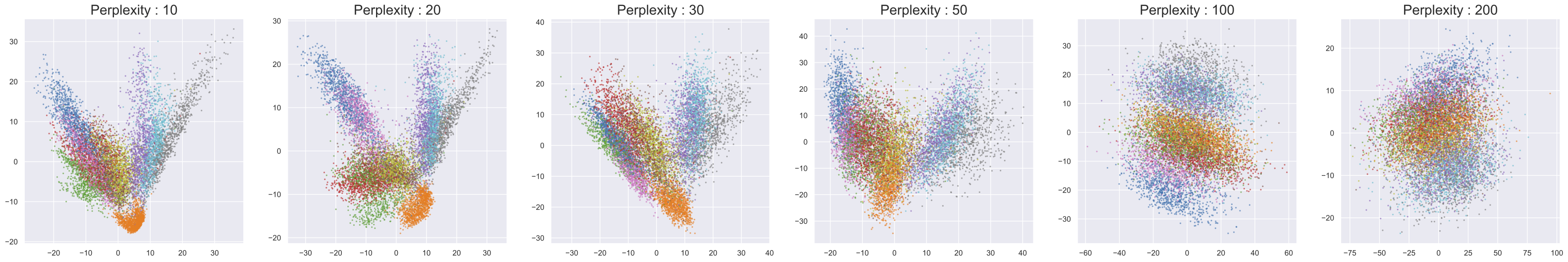}}
\caption{\textit{ccPCA} launched for different values of the perplexity parameter. The latter determines the kernel bandwidths and can be interpreted as the number of effective neighbors of each point \cite{maaten2008tSNE}. As the perplexity grows, the probability of connecting different clusters of digit by sampling through the graph posterior $\mathbb{P}_{\scaleto{\mathcal{P}_X}{4pt}}(\cdot;\bm{K}_{\scaleto{X}{4pt}})$ increases. Therefore clusters are less and less identifiable as the perplexity increases.}
\label{fig:ccPCA_perp}
\end{center}
\end{figure*}

Recall that the \textit{ccPCA} algorithm retrieves the same latent graph as neighbor embedding methods. As shown in \cref{sec:retrieving_DR_methods}, these graphs' distributions depend on the type of prior considered, and take simple forms as follows, when $\bm{\pi}_{\scaleto{X}{4pt}} = \bm{1}$ :
\begin{itemize}
    \item if $\mathcal{P} = B$, $\forall (i,j) \in [n]^2, \: W_{ij} \stackrel{\perp\!\!\!\!\perp}{\sim} \mathcal{B}\left(K_{\scaleto{X}{4pt},ij}/(1 + K_{\scaleto{X}{4pt},ij}) \right)$
    \item if $\mathcal{P} = D$, $\forall i \in [n], \: \W_{i} \stackrel{\perp\!\!\!\!\perp}{\sim} \mathcal{M}\left(1, \bm{K}_{\scaleto{X}{4pt},i}/K_{\scaleto{X}{4pt},i+} \right)$
    \item if $\mathcal{P} = E$, $\W \sim \mathcal{M}\left(n, \bm{K}_{\scaleto{X}{4pt}}/ K_{\scaleto{X}{4pt},++} \right)$
\end{itemize}
and $\bm{K}_{\scaleto{X}{4pt}}$ is the kernel matrix evaluated on the data such that:
$$\forall (i,j) \in [n]^2, \quad \bm{K}_{\scaleto{X}{4pt},ij} = k((\X_{i} - \X_{j})/\tau_{i})$$
where $\bm{\tau} \in \mathbb{R}^n$ is set using an heuristic depending on the method considered \cite{maaten2008tSNE, mcinnes2018umap, tang2016visualizing}. In \cref{fig:ccPCA_perp}, we focus on the effect of the kernel bandwidths on \textit{ccPCA}, choosing the example of t-SNE.

From \cref{fig:ccPCA_perp}, one can notably notice that using a high perplexity leads to a more connected graph and therefore a PCA-like embedding with less degeneracy and no clustering effect. Recall that \textit{ccPCA} computes the same clusters as t-SNE through the CCs of the latent MRF and manage to position t-SNE clusters by focusing on their relative positions (that are filtered by t-SNE). In the case of a connected graph (high perplexity), \textit{ccPCA} will show little advantage over classical PCA since there will not be any cluster to position. Note that this discussion can be extended to other neighbor embedding methods equivalently. Therefore, our probabilistic framework allows us to indentify which part of information is filtered by the posterior graphs with given kernel bandwidth.

\subsection{Quantitative Evaluation of \textit{ccPCA}}\label{sec:quantitative_evaluation}

For quantitative assessment of \textit{ccPCA}, we focused on t-SNE \cite{maaten2008tSNE} and UMAP \cite{mcinnes2018umap} which are the most popular neighbor embedding methods. Note that for these algorithms the initialization is crucial for the global structure of the embeddings as shown in \cite{kobak2021initialization}. In addition to MNIST \cite{deng2012mnist}, we considered the datasets cifar-10, cifar-100 \cite{krizhevsky2009cifar}, fashion-MNIST \cite{xiao2017fashion} as well as the CD8+ T lymphocytes single cell RNA-seq dataset from \cite{kurd2020early}.

We used the quantitative criterion of \cite{lee2015multi} to assess the quality of the embeddings. As mentionned in this paper, the use of this criterion appears as the general consensus in dimension reduction, a field in which building meaningful criteria is tedious. The criterion measures the rescaled average agreement between the K-ary neighbourhoods in the input and output spaces. It is constructed as follows.

We first define the following quantities for $(i,j) \in [n]^2$,
$\rho_{ij} = \mid \{ k: ||\X_{i} - \X_k ||^2_2 < ||\X_{i} - \X_j ||^2_2 \} \mid, r_{ij} = \mid \{ k: ||\Z_{i} - \Z_k ||^2_2 < ||\Z_{i} - \Z_j ||^2_2 \} \mid, \nu_i^K = \{ j : 1 \leq \rho_{ij} \leq K \} \text{ and } \gamma_i^K = \{ j : 1 \leq r_{ij} \leq K \}$. The average K-ary neighbourhood preservation is rescaled to indicate the improvement over a random embedding such that:
\begin{align}\label{def_R}
R_{n}(K) = \frac{(n-1) Q_{n}(K) - K }{n-1-K}
\end{align}
where $Q_{n}(K) = \frac{1}{Kn} \sum_{i=1}^n \mid \nu_i^K \cap \gamma_i^K \mid$, $n$ is the number of data points and $K$ is the hyperparameter that adjusts the scale at which we are looking.

To focus on large-scale structure, K was chosen as either n/4 or n/2. As summarized by \cite{kobak2021initialization}, current practice consists in using PCA or Laplacian eigenmaps as initialization for these algorithms, thus we compare to these strategies. Results are displayed in \cref{tSNE_quantitative_results} and \cref{UMAP_quantitative_results}, each entry being an average over 5 random seeds, with standard deviation displayed below each entry. Note that when not specified, tSNE and UMAP are initialized with an isotropic Gaussian variable. 

These results show that using \textit{ccPCA} is a reliable alternative to PCA and Laplacian eigenmaps for reproducing large-scale neighborhoods.

\vspace{2cm}

\begin{table}[h]
\caption{$100 \times R_{n}(K)$ (\ref{def_R}) for embeddings produced using t-SNE with various initializations.}
\begin{center}
\begin{tabular}{lcccccccc}
\toprule
& \multicolumn{2}{c}{tSNE} & \multicolumn{2}{c}{PCA + tSNE} & \multicolumn{2}{c}{LE + tSNE} & \multicolumn{2}{c}{ccPCA + tSNE} \\
\cmidrule(lr){2-3}\cmidrule(lr){4-5}\cmidrule(lr){6-7}\cmidrule(lr){8-9}
K of K-ary & {n/4} & {n/2}  & {n/4} & {n/2} & {n/4} & {n/2} & {n/4} & {n/2} \\
\midrule
MNIST & 18.7 & 7.4 & 28.4 & 21.9 & 26.7 & 18.5 & $\bm{31.3}$ & $\bm{28.5}$ \\
& $\pm 2.2$ & $\pm 5.1$ & $\pm 0.3$ & $\pm 0.2$ & $\pm 0.7$ & $\pm 0.4$ & $\pm 0.4$ & $\pm 1.2$ 
\vspace{0.1cm}
\\
cifar-10 & 20.3 & 16.4 & $\bm{36.9}$ & 41.9 & 25.8 & 24.1 & 36.4 & $\bm{43.4}$ \\
& $\pm 3.2$ & $\pm 4.8$ & $\pm 0.6$ & $\pm 1.1$ & $\pm 0.6$ & $\pm 1.5$ & $\pm 0.4$ & $\pm 1.6$
\vspace{0.1cm}
\\
cifar-100 & 21.6 & 18.2 & 38.1 & $\bm{47.5}$ & 23.3 & 26.5 & $\bm{39.6}$ & 43.6 \\
& $\pm 3.6$ & $\pm 5.5$ & $\pm 0.4$ & $\pm 0.4$ & $\pm 1.5$ & $\pm 1.8$ & $\pm 0.7$ & $\pm 1.1$
\vspace{0.1cm}
\\
fashion-MNIST & 27.2 & 12.3 & 36.9 & 28.5 & 32.0 & 25.1 & $\bm{41.6}$ & $\bm{35.7}$ \\
& $\pm 4.3$ & $\pm 7.8$ & $\pm 0.1$ & $\pm 0.2$ & $\pm 0.8$ & $\pm 2.2$ & $\pm 0.9$ & $\pm 1.5$
\vspace{0.1cm}
\\
Single Cell data & 25.7 & 22.4 & 37.7 & 29.0 & 28.1 & 31.5 & $\bm{40.1}$ & $\bm{34.6}$ \\
& $\pm 4.8$ & $\pm 10.6$ & $\pm 2.7$ & $\pm 4.7$ & $\pm 1.5$ & $\pm 1.4$ & $\pm 1.7$ & $\pm 2.6$
\\
\bottomrule
\end{tabular}
\end{center}
\label{tSNE_quantitative_results}
\end{table}

\vspace{2cm}

\begin{table}[h]
\caption{$100 \times R_{n}(K)$ (\ref{def_R}) for embeddings produced using UMAP with various initializations.}
\begin{center}
\begin{tabular}{lcccccccc}
\toprule
& \multicolumn{2}{c}{UMAP} & \multicolumn{2}{c}{PCA + UMAP} & \multicolumn{2}{c}{LE + UMAP} & \multicolumn{2}{c}{ccPCA + UMAP} \\
\cmidrule(lr){2-3}\cmidrule(lr){4-5}\cmidrule(lr){6-7}\cmidrule(lr){8-9}
K of K-ary & {n/4} & {n/2}  & {n/4} & {n/2} & {n/4} & {n/2} & {n/4} & {n/2} \\
\midrule
MNIST & 29.5 & 22.7 & $\bm{36.6}$ & 31.1 & 34.6 & 24.9 & 33.4 & $\bm{32.3}$ \\
& $\pm 1.4$ & $\pm 2.2$ & $\pm 0.2$ & $\pm 0.5$ & $\pm 0.2$ & $\pm 0.7$ & $\pm 0.3$ & $\pm 0.5$ 
\vspace{0.1cm}
\\
cifar-10 & 39.2 & 47.6 & 44.3 & $\bm{53.4}$ & 44.2 & 52.6 & $\bm{44.6}$ & 53.2\\
& $\pm 2.6$ & $\pm 1.1$ & $\pm 0.2$ & $\pm 0.1$ & $\pm 0.1$ & $\pm 0.2$ & $\pm 0.2$ & $\pm 0.2$ 
\vspace{0.1cm}
\\
cifar-100 & 41.6 & 42.2 & 45.4 & 45.2 & 44.2 & 43.4 & $\bm{49.9}$ & $\bm{52.9}$\\
& $\pm 1.8$ & $\pm 0.9$ & $\pm 0.2$ & $\pm 0.1$ & $\pm 0.3$ & $\pm 0.1$ & $\pm 0.4$ & $\pm 0.6$ 
\vspace{0.1cm}
\\
fashion-MNIST & 48.7 & 33.6 & 56.2 & 54.3 & 58.1 & 53.4 & $\bm{58.9}$ & $\bm{55.8}$ \\
& $\pm 2.6$ & $\pm 9.5$ & $\pm 0.5$ & $\pm 0.6$ & $\pm 0.5$ & $\pm 0.6$ & $\pm 0.5$ & $\pm 0.3$ 
\vspace{0.1cm}
\\
Single Cell data & 39.5 & 34.3 & 52.3 & 47.2 & $\bm{55.9}$ & 45.7 & 53.6 & $\bm{53.9}$\\
& $\pm 1.4$ & $\pm 6.1$ & $\pm 0.8$ & $\pm 6.9$ & $\pm 0.3$ & $\pm 0.9$ & $\pm 0.3$ & $\pm 1.3$ 
\vspace{0.1cm}
\\
\bottomrule
\end{tabular}
\end{center}
\label{UMAP_quantitative_results}
\end{table}


\begin{thebibliography}{10}

\bibitem{anders2018dissecting}
Friedrich Anders, Cristina Chiappini, Basilio~Xavier Santiago, Gal
  Matijevi{\v{c}}, Anna~B Queiroz, Matthias Steinmetz, and Guillaume Guiglion.
\newblock Dissecting stellar chemical abundance space with t-sne.
\newblock {\em Astronomy \& Astrophysics}, 619:A125, 2018.

\bibitem{anowar2021conceptual}
Farzana Anowar, Samira Sadaoui, and Bassant Selim.
\newblock Conceptual and empirical comparison of dimensionality reduction
  algorithms ({PCA, KPCA, LDA, MDS, SVD, LLE, ISOMAP, LE, ICA, t-SNE)}.
\newblock {\em Computer Science Review}, 40:100378, 2021.

\bibitem{arora2018analysis}
Sanjeev Arora, Wei Hu, and Pravesh~K Kothari.
\newblock An analysis of the t-sne algorithm for data visualization.
\newblock In {\em Conference On Learning Theory}, pages 1455--1462. PMLR, 2018.

\bibitem{balasubramanian2002isomap}
Mukund Balasubramanian, Eric~L Schwartz, Joshua~B Tenenbaum, Vin de~Silva, and
  John~C Langford.
\newblock The isomap algorithm and topological stability.
\newblock {\em Science}, 295(5552):7--7, 2002.

\bibitem{belkin2003laplacian}
Mikhail Belkin and Partha Niyogi.
\newblock Laplacian eigenmaps for dimensionality reduction and data
  representation.
\newblock {\em Neural computation}, 15(6):1373--1396, 2003.

\bibitem{bohm2020unifying}
Jan~Niklas B{\"o}hm, Philipp Berens, and Dmitry Kobak.
\newblock A unifying perspective on neighbor embeddings along the
  attraction-repulsion spectrum.
\newblock {\em arXiv preprint arXiv:2007.08902}, 2020.

\bibitem{carreira2010elastic}
Miguel~A Carreira-Perpin{\'a}n.
\newblock The elastic embedding algorithm for dimensionality reduction.
\newblock In {\em ICML}, volume~10, pages 167--174. Citeseer, 2010.

\bibitem{chan2018t}
David~M Chan, Roshan Rao, Forrest Huang, and John~F Canny.
\newblock t-sne-cuda: Gpu-accelerated t-sne and its applications to modern
  data.
\newblock In {\em 2018 30th International Symposium on Computer Architecture
  and High Performance Computing (SBAC-PAD)}, pages 330--338. IEEE, 2018.

\bibitem{Chung97}
F.~R.~K. Chung.
\newblock {\em Spectral Graph Theory}.
\newblock American Mathematical Society, 1997.

\bibitem{clifford1990markov}
Peter Clifford.
\newblock Markov random fields in statistics.
\newblock {\em Disorder in physical systems: A volume in honour of John M.
  Hammersley}, pages 19--32, 1990.

\bibitem{coenen2019understanding}
Andy Coenen and Adam Pearce.
\newblock Understanding umap.
\newblock {\em Google PAIR}, 2019.

\bibitem{coifman2006diffusion}
Ronald~R Coifman and St{\'e}phane Lafon.
\newblock Diffusion maps.
\newblock {\em Applied and computational harmonic analysis}, 21(1):5--30, 2006.

\bibitem{csardi2006igraph}
Gabor Csardi, Tamas Nepusz, et~al.
\newblock The igraph software package for complex network research.
\newblock {\em InterJournal, complex systems}, 1695(5):1--9, 2006.

\bibitem{deng2012mnist}
Li~Deng.
\newblock The mnist database of handwritten digit images for machine learning
  research.
\newblock {\em IEEE Signal Processing Magazine}, 29(6):141--142, 2012.

\bibitem{donoho2000high}
David~L Donoho.
\newblock High-dimensional data analysis: The curses and blessings of
  dimensionality.
\newblock {\em AMS math challenges lecture}, 1(2000):32, 2000.

\bibitem{ham2004kernel}
Jihun Ham, Daniel~D Lee, Sebastian Mika, and Bernhard Sch{\"o}lkopf.
\newblock A kernel view of the dimensionality reduction of manifolds.
\newblock In {\em Proceedings of the twenty-first international conference on
  Machine learning}, page~47, 2004.

\bibitem{hinton2002stochastic}
Geoffrey Hinton and Sam~T Roweis.
\newblock Stochastic neighbor embedding.
\newblock In {\em NIPS}, volume~15, pages 833--840. Citeseer, 2002.

\bibitem{NIPS2002SNE}
Geoffrey~E Hinton and Sam~T. Roweis.
\newblock Stochastic neighbor embedding.
\newblock In S.~Becker, S.~Thrun, and K.~Obermayer, editors, {\em Advances in
  Neural Information Processing Systems 15}, pages 857--864. MIT Press, 2003.

\bibitem{kobak2019art}
Dmitry Kobak and Philipp Berens.
\newblock The art of using t-sne for single-cell transcriptomics.
\newblock {\em Nature communications}, 10(1):1--14, 2019.

\bibitem{kobak2019heavy}
Dmitry Kobak, George Linderman, Stefan Steinerberger, Yuval Kluger, and Philipp
  Berens.
\newblock Heavy-tailed kernels reveal a finer cluster structure in t-sne
  visualisations.
\newblock In {\em Joint European Conference on Machine Learning and Knowledge
  Discovery in Databases}, pages 124--139. Springer, 2019.

\bibitem{kobak2021initialization}
Dmitry Kobak and George~C Linderman.
\newblock Initialization is critical for preserving global data structure in
  both t-sne and umap.
\newblock {\em Nature biotechnology}, 39(2):156--157, 2021.

\bibitem{krizhevsky2009cifar}
Alex Krizhevsky, Vinod Nair, and Geoffrey Hinton.
\newblock Cifar-10 and cifar-100 datasets.
\newblock {\em URl: https://www. cs. toronto. edu/kriz/cifar. html}, 6(1):1,
  2009.

\bibitem{kruskal1978multidimensional}
Joseph~B Kruskal.
\newblock {\em Multidimensional scaling}.
\newblock Number~11. Sage, 1978.

\bibitem{kurd2020early}
Nadia~S Kurd, Zhaoren He, Tiani~L Louis, J~Justin Milner, Kyla~D Omilusik,
  Wenhao Jin, Matthew~S Tsai, Christella~E Widjaja, Jad~N Kanbar, Jocelyn~G
  Olvera, et~al.
\newblock Early precursors and molecular determinants of tissue-resident memory
  cd8+ t lymphocytes revealed by single-cell rna sequencing.
\newblock {\em Science immunology}, 5(47):eaaz6894, 2020.

\bibitem{lee2015multi}
John~A Lee, Diego~H Peluffo-Ord{\'o}{\~n}ez, and Michel Verleysen.
\newblock Multi-scale similarities in stochastic neighbour embedding: Reducing
  dimensionality while preserving both local and global structure.
\newblock {\em Neurocomputing}, 169:246--261, 2015.

\bibitem{li2020high}
Tianxi Li, Cheng Qian, Elizaveta Levina, and Ji~Zhu.
\newblock High-dimensional gaussian graphical models on network-linked data.
\newblock {\em J. Mach. Learn. Res.}, 21:74--1, 2020.

\bibitem{li2017application}
Wentian Li, Jane~E Cerise, Yaning Yang, and Henry Han.
\newblock Application of t-sne to human genetic data.
\newblock {\em Journal of bioinformatics and computational biology},
  15(04):1750017, 2017.

\bibitem{linderman2019fast}
George~C Linderman, Manas Rachh, Jeremy~G Hoskins, Stefan Steinerberger, and
  Yuval Kluger.
\newblock Fast interpolation-based t-sne for improved visualization of
  single-cell rna-seq data.
\newblock {\em Nature methods}, 16(3):243--245, 2019.

\bibitem{linderman2019clustering}
George~C Linderman and Stefan Steinerberger.
\newblock Clustering with t-sne, provably.
\newblock {\em SIAM Journal on Mathematics of Data Science}, 1(2):313--332,
  2019.

\bibitem{mcinnes2018umap}
Leland McInnes, John Healy, and James Melville.
\newblock Umap: Uniform manifold approximation and projection for dimension
  reduction.
\newblock {\em arXiv preprint arXiv:1802.03426}, 2018.

\bibitem{pearson1901liii}
Karl Pearson.
\newblock Liii. on lines and planes of closest fit to systems of points in
  space.
\newblock {\em The London, Edinburgh, and Dublin philosophical magazine and
  journal of science}, 2(11):559--572, 1901.

\bibitem{pedregosa2011scikit}
Fabian Pedregosa, Ga{\"e}l Varoquaux, Alexandre Gramfort, Vincent Michel,
  Bertrand Thirion, Olivier Grisel, Mathieu Blondel, Peter Prettenhofer, Ron
  Weiss, Vincent Dubourg, et~al.
\newblock Scikit-learn: Machine learning in python.
\newblock {\em the Journal of machine Learning research}, 12:2825--2830, 2011.

\bibitem{pezzotti2019gpgpu}
Nicola Pezzotti, Julian Thijssen, Alexander Mordvintsev, Thomas H{\"o}llt,
  Baldur Van~Lew, Boudewijn~PF Lelieveldt, Elmar Eisemann, and Anna Vilanova.
\newblock Gpgpu linear complexity t-sne optimization.
\newblock {\em IEEE transactions on visualization and computer graphics},
  26(1):1172--1181, 2019.

\bibitem{polivcar2019opentsne}
Pavlin~G Poli{\v{c}}ar, Martin Stra{\v{z}}ar, and Bla{\v{z}} Zupan.
\newblock opentsne: a modular python library for t-sne dimensionality reduction
  and embedding.
\newblock {\em BioRxiv}, page 731877, 2019.

\bibitem{rue2005gaussian}
Havard Rue and Leonhard Held.
\newblock {\em Gaussian Markov random fields: theory and applications}.
\newblock CRC press, 2005.

\bibitem{tang2016visualizing}
Jian Tang, Jingzhou Liu, Ming Zhang, and Qiaozhu Mei.
\newblock Visualizing large-scale and high-dimensional data.
\newblock In {\em Proceedings of the 25th international conference on world
  wide web}, pages 287--297, 2016.

\bibitem{tipping1999probabilistic}
Michael~E Tipping and Christopher~M Bishop.
\newblock Probabilistic principal component analysis.
\newblock {\em Journal of the Royal Statistical Society: Series B (Statistical
  Methodology)}, 61(3):611--622, 1999.

\bibitem{maaten2008tSNE}
L.J.P. {van der Maaten} and G.E. Hinton.
\newblock Visualizing high-dimensional data using t-sne.
\newblock {\em Journal of Machine Learning Research}, 9(nov):2579--2605, 2008.
\newblock Pagination: 27.

\bibitem{wainwright2008graphical}
Martin~J Wainwright and Michael~Irwin Jordan.
\newblock {\em Graphical models, exponential families, and variational
  inference}.
\newblock Now Publishers Inc, 2008.

\bibitem{wang2021understanding}
Yingfan Wang, Haiyang Huang, Cynthia Rudin, and Yaron Shaposhnik.
\newblock Understanding how dimension reduction tools work: an empirical
  approach to deciphering t-sne, umap, trimap, and pacmap for data
  visualization.
\newblock {\em J Mach. Learn. Res}, 22:1--73, 2021.

\bibitem{wattenberg2016use}
Martin Wattenberg, Fernanda Vi{\'e}gas, and Ian Johnson.
\newblock How to use t-sne effectively.
\newblock {\em Distill}, 1(10):e2, 2016.

\bibitem{xiao2017fashion}
Han Xiao, Kashif Rasul, and Roland Vollgraf.
\newblock Fashion-mnist: a novel image dataset for benchmarking machine
  learning algorithms.
\newblock {\em arXiv preprint arXiv:1708.07747}, 2017.

\end{thebibliography}
\end{document}